\documentclass[a4paper]{amsart}
\usepackage{graphicx,color}
\usepackage[psamsfonts]{amssymb}
\usepackage{hyperref}
\usepackage{enumerate}

\newtheorem{theorem}{Theorem}[section]
\newtheorem{corollary}[theorem]{Corollary}
\newtheorem{lemma}[theorem]{Lemma}
\newtheorem{proposition}[theorem]{Proposition}
\newtheorem{claim}{Claim}[theorem]
\theoremstyle{definition}

\newtheorem{remark}[theorem]{Remark}
\newtheorem{example}[theorem]{Example}

\numberwithin{equation}{section}

\newcommand{\set}[1]{\left\{#1\right\}}

\newcommand{\inte}{{\mathop{\mathrm{int}\,}}}
\newcommand{\cl}{{\mathop{\mathrm{cl}\,}}}
\newcommand{\relbd}{{\mathop{\mathrm{relbd}\,}}}
\newcommand{\relint}{{\mathop{\mathrm{relint}\,}}}

\newcommand{\conv}{{\mathop{\mathrm{conv}\,}}}
\newcommand{\cone}{{\mathop{\mathrm{cone}\,}}}

\newcommand{\supp}{{\mathop{\mathrm{supp}\,}}}

\newcommand{\len}{{\mathop{\lambda_1}}}
\newcommand{\area}{{\mathop{\lambda_2}}}
\newcommand{\voln}{{\mathop{\lambda_n}}}

\newcommand{\pa}{{\partial}}

\newcommand{\Real}{\mathbb{R}}

\newcommand{\ee}{{\varepsilon}}
\newcommand{\al}{{\alpha}}
\newcommand{\be}{{\beta}}
\newcommand{\de}{{\delta}}
\newcommand{\ga}{{\gamma}}

\newcommand{\la}{{\lambda}}
\newcommand{\te}{{\theta}}
\newcommand{\Om}{{\Omega}}
\newcommand{\Si}{{\Sigma}}
\newcommand{\De}{{\Delta}}
\newcommand{\calA}{{\mathcal A}}

\newcommand{\F}{{\mathcal A}}
\newcommand{\G}{{\mathcal B}}
\newcommand{\cK}{{\mathcal K}}
\newcommand{\cL}{{\mathcal L}}
\newcommand{\cT}{{\mathcal T}}

\newcommand{\jcov}[1]{g_{#1}}

\newcommand{\ang}[2]{{\{#1\leq\theta\leq#2\}}}

\begin{document}
\title[The cross covariogram of pairs of polygons]{The cross covariogram of a pair of polygons determines both polygons, with a few exceptions}
\author{Gabriele Bianchi}
\address{Dipartimento di Matematica, Universit\`a di Firenze, 
Viale Morgagni 67/A, Firenze, Italy I-50134}
\email{gabriele.bianchi@unifi.it}
\date{\today}

\begin{abstract}
The cross covariogram $g_{K,L}$ of two convex sets
$K$ and $L$ in $\Real^n$ is the function which associates to each $x\in
\Real^n$ the volume of $K\cap(L+x)$.
Very recently Averkov and Bianchi~\cite{AB2} have confirmed Matheron's conjecture on the covariogram problem, that asserts that any planar convex body $K$ is determined by the knowledge of $g_{K,K}$. 
The problem of determining the sets from their covariogram is relevant in probability, in statistical shape recognition and in the determination of the atomic structure of a quasicrystal from X-ray diffraction images.
We prove that when  $K$ and $L$ are convex polygons (and also when $K$ and $L$ are planar convex cones)  $g_{K,L}$ determines
both $K$ and $L$, up to a described family of exceptions. These  results  imply that, when $K$ and $L$ are in these classes, the information provided by the cross covariogram is so rich as to  determine not only one unknown body, as required by Matheron's conjecture, but two bodies, with a few classified exceptions.  These results are also used by Bianchi~\cite{B3d} to prove that any convex polytope $P$ in $\Real^3$ is determined by $g_{P,P}$.
\end{abstract}
\maketitle

\section{Introduction}
Let $K$ and $L$ be convex sets in  $\Real^n$, $n \geq 2$, and let $\voln$ stand for the $n$-dimensional Lebesgue measure. 
The \emph{cross covariogram} $\jcov{K,L}$ of $K$ and $L$ is the
function defined by
\[
\jcov{K,L}(x)=\voln(K\cap (L+x)),
\]
where $x\in\Real^n$ is such that $\voln(K\cap (L+x))$ is finite. This function, introduced by Cabo and Janssen~\cite{CJ}, coincides with the convolution of  the characteristic function $\mathbf{1}_K$ of $K$ with the characteristic function $\mathbf{1}_{-L}$ of the reflection of $L$ in the origin, that is,
\begin{equation}\label{convoluzione}
g_{K,L} =\mathbf{1}_K\ast \mathbf{1}_{-L}.
\end{equation}

The function $g_{K,K}$  was introduced by G.~Matheron in his book~\cite[Section~4.3]{M} on random sets, is denoted by $g_K$ and is called \emph{covariogram} or \emph{set covariance} of $K$. 
The covariogram $g_K$ is clearly unchanged by a translation or a reflection 
of $K$. (The term \emph{reflection} will always mean reflection 
in a point.) A \emph{convex body} in $\Real^n$ is a convex compact set with non-empty interior. In 1986 Matheron~\cite[p.~20]{M2}  asked the following question and conjectured a positive answer for the case $n=2$.
\smallskip

\textbf{Covariogram problem.} \emph{Does $g_K$ determine a
convex body $K$, among all convex bodies, up to translations and
reflections?}
\smallskip

Even for $n=2$ this conjecture has been completely settled only very recently, by  Averkov and Bianchi~\cite{AB2}. 
It is known that the covariogram problem is equivalent to any of the following problems (see~\cite{AB2} for a detailed explanation).
\begin{enumerate}[{P}1.]
\item\label{chord_lengths} Determine a convex  body $K$ by the knowledge, for each unit vector $u$ in $\Real^n$, of the distribution of the lengths of the chords of $K$ parallel to $u$.
\item\label{distribution_XY} Determine a convex  body $K$ by the  distribution of $X-Y$, where $X$ and $Y$ are independent random variables uniformly distributed over $K$.
\item\label{phase_retrieval} Determine  the characteristic function $\mathbf{1}_K$ of a convex body $K$ from  the modulus of its Fourier transform $\widehat{\mathbf{1}_K}$.
\end{enumerate}
Chord-length distributions  are of wide interest beyond mathematics and are common in \emph{stereology, statistical shape recognition and image analysis,} where properties of an unknown body
have to be inferred from chord length measurements; see Schmitt~\cite{S},
Cabo and Baddeley~\cite{CB}, Serra~\cite{Se} and Mazzolo, Roesslinger and Gille~\cite{MRG03}. 
Problem P\ref{distribution_XY} was asked by Adler and Pyke~\cite{AP1} in 1991. The same authors~\cite{AP2} find the
covariogram problem relevant also in the study of scanning Brownian
processes and of the equivalence of measures induced by these
processes for different base sets.
Problem P\ref{phase_retrieval} is a special case of the \emph{phase retrieval problem}, where $\mathbf{1}_K$ is replaced by a function with compact support. The phase retrieval problem has applications in \emph{X-ray crystallography, optics, electron microscopy} and other areas, references to which may be found in~\cite{BSV}. 
Very recently, Baake and Grimm~\cite{BaakeGrimm} have proved that the covariogram problem is particularly relevant for finding the atomic structure of a \emph{quasicrystal} from its $X$-ray diffraction image. When a quasicrystal  fits into the  so-called \emph{cut-and-project scheme}, the determination of its atomic structure $S$ requires the knowledge of an unknown ``window'' $W$, which in many important cases is a convex body. The covariogram problem enters at this point, since $g_W$ can be obtained from the diffraction image of $S$. 

We have already mentioned that the covariogram problem has a positive answer in the plane. In higher dimensions a complete answer is known only when $K$ is a  convex polytope. Bianchi~\cite{B2}  proved that in $\Real^n$, for every $n\geq 4$, there are pairs of convex polytopes with equal covariograms  which are not translations or reflections of each other. On the other hand, the answer to the covariogram problem for a three-dimensional convex polytope is positive, as proved by  Bianchi~\cite{B3d}. 

Cabo and Janssen~\cite{CJ} prove that $g_{C,-C}$ determines every regular (equal to the closure of its interior) compact set $C$ in $\Real^n$, $n\geq2$. This result clearly implies that the covariogram determines each centrally symmetric regular compact set in $\Real^n$.  In general, the convexity assumption in the covariogram problem is needed, since there exist examples of  non-convex polyominoes which are neither translations nor reflections of each other and have equal covariograms; see Gardner, Gronchi and Zong~\cite{GGZ}. However, a planar convex body is determined by its covariogram in a class that is much larger than that of convex bodies; see Benassi, Bianchi and D'Ercole~\cite{BenassiBianchiDErcole}.

When $K$ is a planar convex body, the information provided by $g_K$ seems to be  richer than  is necessary to determine $K$. For instance, for a convex body $K$ whose boundary is  $C^2$ regular and has non-zero curvature, Averkov and Bianchi~\cite{AB} indicate some subsets of the 
support of $g_K$, with arbitrarily small Lebesgue measure, such that $g_K$ restricted to those subsets identifies $K$. 
In this paper we investigate this richness of information from a different point of view, trying to understand which information $g_{K,L}$ carries about the \emph{two} convex sets $K$ and $L$. We are able to provide a complete answer when $K$ and $L$ are convex polygons, and also when they are planar convex cones. The obtained results  imply that the information provided by the cross covariogram, when $K$ and $L$ are in these classes, is so rich as to  determine not only one unknown body, as required by Matheron's conjecture, but two bodies, with a few classified exceptions. 

In order to generalise the covariogram problem to the case of the cross covariogram, we  observe that the translation of $K$ and $L$ by the same vector, and the substitution of $K$ with $-L$ and of $L$ with $-K$, leave $g_{K,L}$ unchanged.  
Let $K$, $L$, $K'$ and $L'$ be  convex sets in $\Real^2$. 
We call $(K,L)$ and $(K',L')$ \emph{trivial associates} when one pair is obtained by the other one via a combination of the two operations above, that is, when either $(K,L)=(K'+x,L'+x)$ or $(K,L)=(-L'+x,-K'+x)$, for some $x\in \Real^2$.

\textbf{Cross covariogram problem.} 
\emph{Does $g_{K,L}$ determine the pair of closed convex sets $(K,L)$, among all pairs of closed convex sets, up to trivial associates?}

Assume $K$ and $L$ convex polygons. In this case the answer to the cross covariogram problem is negative as Examples~\ref{parall} and~\ref{parall_due} show (see Figures~\ref{fig_parall} and~\ref{fig_parall_due}). For each choice of some real parameters there exist four  pairs of parallelograms $(\cK_1,\cL_1),\dots,(\cK_4,\cL_4)$ such that, for $i=1,3$, $g_{\cK_i,\cL_i}=g_{\cK_{i+1},\cL_{i+1}}$ but $(\cK_i,\cL_i)$ is not a trivial associate of $(\cK_{i+1},\cL_{i+1})$. 
Theorem~\ref{cov_congiunto_poligoni} proves that, up to an affine
transformation, the previous counterexamples are the only ones.

\begin{theorem}\label{cov_congiunto_poligoni}
Let $K$ and $L$ be convex polygons and $K'$ and $L'$ be planar closed convex sets
with $g_{K,L}=g_{K',L'}$. Assume that there exist no affine transformation $\cT$ and no different  indices $i,j$, with either $i,j\in\{1,2\}$ or  $i,j\in\{3,4\}$, such that $(\cT K,\cT L)$ and $(\cT K',\cT L')$ are trivial associates of $(\cK_i,\cL_i)$ and of $(\cK_j,\cL_j)$, respectively.
Then  $(K,L)$ is a trivial associate of $(K',L')$.
\end{theorem}

Theorem~\ref{cov_congiunto_poligoni} has a probabilistic interpretation in terms of a generalisation of Problem~P\ref{distribution_XY}. 
It implies that the distribution of the difference $X-Y$ of  two independent random
variables $X$ and $Y$, with $X$ uniformly distributed over a
convex polygon $K$ and $Y$ uniformly distributed over a convex polygon $L$, together with $\area(K)\,\area(L)$, 
determines both $K$ and $L$, up to some inherent ambiguities, with a few exceptions. This result is a consequence of  Theorem~\ref{cov_congiunto_poligoni} because  the probability distribution of $X-Y$ is $g_{K,L}/(\area(K)\,\area(L))$, by \eqref{convoluzione}. 

It would be interesting to understand if  a result similar to Theorem~\ref{cov_congiunto_poligoni} holds for other classes of planar convex bodies, for instance for  bodies with $C^2$ boundary.

Another motivation for Theorem~\ref{cov_congiunto_poligoni} comes from the proof of the positive answer to the covariogram problem for  three-dimensional convex polytopes mentioned above, of which Theorem~\ref{cov_congiunto_poligoni} constitutes a crucial step. The two problems are connected because when $K$ and $L$ are parallel antipodal facets of  a three-dimensional convex polytope $P$,  $g_P$ provides $g_{K,L}$, and Theorem~\ref{cov_congiunto_poligoni} helps to determine those pairs of facets; see \cite{B3d} for the details. 

The previous theorem has something to say also regarding the symmetries of $g_{K,L}$. Let $F$ and $G$ be convex bodies in $\Real^n$. It is evident that $g_F(x)=g_F(-x)$ for every $x\in\Real^n$, but the cross covariogram is not always an even function. (In general, one only has $g_{F,G}(-x)=g_{G,F}(x)=g_{-F,-G}(x)$.)

\begin{corollary}\label{symmetry_crosscov}
Let $K$ and $L$ be convex polygons. Then there exists $z\in\Real^2$ such that $g_{K,L}(z+x)=g_{K,L}(z-x)$ for every $x\in\Real^2$ if and only if either $K=L+z$ or both $K$ and $L+z$ are centrally symmetric with equal center.
\end{corollary}

Mani-Levitska~\cite{ML} saw the study of the cross covariogram problem for pairs of polyhedral convex cones in $\Real^3$ as a step towards the solution of the covariogram problem for polytopes in $\Real^3$, and indeed~\cite{B3d} contains some  results in  this direction. These results for cones in $\Real^3$ are not exhaustive, but  the cross covariogram problem for planar convex cones can be completely understood.
Let $A$, $A'$, $B$ and $B'$ be convex cones  in $\Real^2$ with apex the origin $O$. Assume $\inte A\cap \inte B=\emptyset$, because otherwise $g_{A,B}$ is nowhere finite.
Since the cones have apex  $O$,  $(A,B)$ and $(A',B')$ are trivial associates if and only if $\set{A,-B}=\set{A',-B'}$.
Example~\ref{example_cones} (see Figure~\ref{fig_settori_angolari}) presents two different pairs of convex cones $(\F_1,\G_1)$ and $(\F_2,\G_2)$ with equal cross covariogram  and Theorem~\ref{cov_cong_settori_angolari} proves that, up to affine transformations, these are the only counterexamples.

\begin{theorem}\label{cov_cong_settori_angolari}
Let $A$, $B$, $A'$ and $B'$ be pointed closed convex cones in $\Real^2$ with non-empty interior and apex the origin $O$, such that $\inte A\cap \inte B=\emptyset$.
The identity $\jcov{A,B}=\jcov{A',B'}$ holds true if and only if one of the following alternatives occurs:
\begin{enumerate}
\item\label{acaso0} $\set{A,-B}=\set{A',-B'}$;
\item\label{acaso12} there exist a linear transformation $\cT$ and $i,j\in\set{1,2}$, $i\neq j$, such that
\begin{equation}\label{acaso1}
\{\cT A,-\cT B\}=\{\F_i,-\G_i\}\quad\text{and}\quad \{\cT A',-\cT B'\}=\{\F_j,-\G_j\}.
\end{equation}
\end{enumerate}
\end{theorem}

A crucial notion in the proof of  Theorem~\ref{cov_congiunto_poligoni} is that of \emph{synisothetic} pairs of polytopes, introduced by~\cite{B3d} and explained in Section~\ref{notations}.
In Section~\ref{sec_crosscov_synisothesis} we prove that, up to affine transformations, $(\cK_1,\cL_1)$ and $(\cK_2,\cL_2)$ are the only pairs of convex polygons with equal cross covariogram which are not synisothetic.
To establish this result we use also Theorem~\ref{cov_cong_settori_angolari}, whose proof is contained in Section~\ref{sec_coni}. 
The proofs of Theorem~\ref{cov_congiunto_poligoni} and of Corollary~\ref{symmetry_crosscov} are contained in Section~\ref{sec_teorema_poligoni}. 
We conclude by mentioning that Lemma~\ref{coni_finale} is a technical result which may be of interest by itself.


\section{Definitions, notations and preliminaries}\label{notations}
As usual, $S^{n-1}$ denotes the unit sphere in $\Real^n$, centred at 
the origin $O$.  For $x$, $y\in \Real^n$,  $\|x\|$ is the Euclidean norm of $x$ and $x\cdot y$ denotes scalar product. For $\de>0$, $B(x,\de)$ denotes the open ball in $\Real^n$ centred at $x$ and with radius $\de$. 
For $\te\in[0,2\pi]$ we write $u(\te)=(\cos\te,\sin\te)\in S^1$.

If $A\subset \Real^n$  we denote by $\inte A$, $\cl A$, 
$\partial A$ and
$\conv A$ the \emph{interior}, \emph{closure}, \emph{boundary} and
\emph{convex hull} of $A$, respectively. 
The \emph{symmetric difference} of the sets $A$ and $B$ is defined by $A \De B=(A\setminus B)\cup(B\setminus A)$. The \emph{Minkowski sum} of $A$ and $B$ is
\[
 A+B=\{x+y : x\in A, y\in B\}.
\]
We write $\la_k$ for $k$-dimensional Lebesgue measure in $\Real^n$, where $k = 1,\dots,n$, and where we identify $\la_k$ with $k$-dimensional Hausdorff measure.

A {\it convex body} $K\subset\Real^n$ is a compact convex set with non-empty interior.
The symbols $\relbd K$ and $\relint K$ indicate respectively the \emph{relative boundary} and  the \emph{relative interior} of $K$. 
The \emph{difference body} of $K$ is defined by $DK=K+(-K)$.
The \emph{support function} of $K$ is defined, for $x\in\Real^n$, by 
$
h_K(x)=\sup \set{x\cdot y : y\in K}.
$
Given $x,y\in\Real^n$, we write $[x,y]$ for the line segment with endpoints $x$ and $y$. When $K$ is a planar convex body and $a,b\in\pa K$, the symbol ${(a,b)}_{\pa K}$  denotes the set of points in $\pa K$ which strictly follow $a$ and strictly precede $b$  in counterclockwise order on $\pa K$, and ${[a,b]}_{\pa K}$ denotes ${(a,b)}_{\pa K}\cup\{a,b\}$. Given an arc $\Om\subset\pa K$ with $\cl \Om={[a,b]}_{\pa K}$, we  call 
$a$ the \emph{lower endpoint} of $\Om$,  $b$  its \emph{upper endpoint}, and, with a small abuse of notation,  we will call $(a,b)_{\pa K}$  the \emph{relative interior} of $\Om$. Given $u,v\in S^1$, $v\geq u$ means $v\in(u,-u)_{S^1}\cup\{u\}$, while $v>u$ means $v\in(u,-u)_{S^1}\cup\{-u\}$.

If $F$ is a face of a convex polytope $P$ in $\Real^n$, the \emph{normal cone} of $P$ at $F$ is denoted by
$N(P,F)$ and is the set of all outer normal vectors to $P$ at $x$,
where $x\in\relint F$, together with $O$. The \emph{support cone of $P$ at $F$} is  the set
\[
\cone(P,F)=\set{\mu(y-x)\;:\;y\in P\;,\;\mu\geq0},
\]
where $x\in \relint F$. Neither definitions  depend on the choice of $x$.  
If $u\in S^{n-1}$, the \emph{exposed face of $P$ in direction $u$} is 
\[
P_u=\{x\in P : x\cdot u=h_P(u)\}.
\]
It is the unique proper face of $P$ such that the relative interior of its normal cone contains $u$. \cite[Th.~1.7.5(c)]{Sc} proves that, when $K$ and $L$ are convex polygons and $u\in S^1$, we have 
\begin{equation}\label{facce_corpodiff}
(K+(-L))_u=K_u+(-L)_u.
\end{equation}

In this paper the term \emph{cone} always means cone with apex $O$. 
A convex cone is \emph{pointed} if its apex is a vertex. Let  $(\rho,\te)$ denote polar coordinates in $\Real^2$ and let $\al,\be\in[0,2\pi]$, with $\al<\be$. For brevity we write $\set{\te=\al}$ for the ray $\set{(\rho,\te):\,\te=\al}$ 
and $\ang{\al}{\be}$ for the  cone $\{(\rho,\te):\,\te\in[\alpha,\beta]\}$. 

Let $\supp f$ denote the support of the function $f$. Let $K$ and $L$ be convex polygons and let $A$ and $B$ be closed convex cones in $\Real^2$ with $\inte A\cap \inte B=\emptyset$. It is easy to prove that
\begin{equation}\label{support}
\supp g_{K,L}=K+(-L) \quad\text{and}\quad\supp g_{A,B}=A+(-B)=\conv(A\cup(-B)).
\end{equation}
It can be proved, by using the Minkowski inequality as in \cite[p.~410-411]{Sc}, that ${g_{K,L}}^{{1}/{2}}$ is concave on its support.

\textit{Synisothesis.}
Let $P$ and $Q$ be convex polytopes in $\Real^n$, let $F$ be a proper face of $P$, and let $G$ be a proper face of $Q$. We say that $F$ and $G$ are {\em isothetic} if $G$ is a translate of $F$ and
\[
\cone(P,F)=\cone(Q,G).
\]
Given convex polytopes $P_1$, $P_2$, $Q_1$ and $Q_2$ in $\Real^n$ we
say that $(P_1,P_2)$ and $(Q_1,Q_2)$ are \emph{synisothetic} if
given any proper face $F$ of $P_1$ or of $P_2$ there is a
proper face $G$  of $Q_1$ or of $Q_2$ (and conversely) such
that $F$ and $G$ are isothetic.

The notion of synisothesis will play a central role in this paper. Let $K$, $K'$, $L$ and $L'$ be convex polygons. 
It is convenient for later use to express the  synisothesis of $(K,-L)$ and $(K',-L')$  in two equivalent ways. 
It is clear that $(K,-L)$ and $(K',-L')$ are synisothetic if and only if, for each $u\in S^1$, one of the following properties hold:
\begin{align} 
&\text{$K_u$ is isothetic to  $K'_u$  and  
$(-L)_u$ is isothetic to $(-L')_u$; }\label{alt_uno} \\ 
&\text{$K_u$ is isothetic to  $(-L')_u$  and  
$(-L)_u$ is isothetic to $K'_u$.}\label{alt_due}
\end{align}
Let $u\in S^1$ and $F$ and $F'$ be convex polygons. Observe that  $F_u$ is isothetic to $F'_u$ if and only if either both are edges  of equal length or both are vertices with equal support cones. Thus,  $(K,-L)$ and $(K',-L')$ are synisothetic if and only if the following equalities hold for each $u\in S^1$: 
\begin{align}
\set{\len(K_u),\len((-L)_u)}
&=\set{\len(K'_u),\len((-L')_u)};\quad\label{sin_poligoni1}\\
 \set{\cone(K,K_u),\cone(-L,(-L)_u)}
&=\set{\cone(K',K'_u),\cone(-L',(-L')_u)}.\label{sin_poligoni2}
\end{align}

\section{The cross covariogram problem for planar convex cones}\label{sec_coni}
Let us introduce the counterexample to the cross covariogram problem for cones.
\begin{example}\label{example_cones}
Let $\F_1=\ang{0}{3\pi/4}$, $\G_1=-\ang{{\pi}/4}{{\pi}/2}$,
$\F_2=\ang{0}{{\pi}/{4}}$ and  $\G_2=-\ang{{\pi}/2}{{3\pi}/4}$; see Fig.~\ref{fig_settori_angolari}.
\begin{figure}
\begin{center}
\begin{picture}(0,0)%
\includegraphics{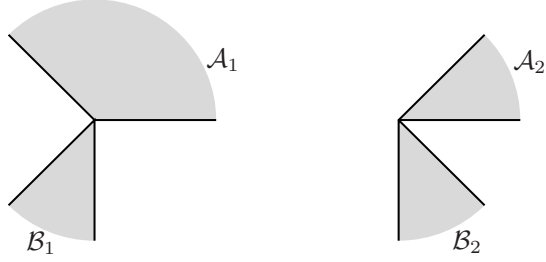}%
\end{picture}%
\setlength{\unitlength}{3315sp}%
\begingroup\makeatletter\ifx\SetFigFontNFSS\undefined%
\gdef\SetFigFontNFSS#1#2#3#4#5{%
  \reset@font\fontsize{#1}{#2pt}%
  \fontfamily{#3}\fontseries{#4}\fontshape{#5}%
  \selectfont}%
\fi\endgroup%
\begin{picture}(3830,1940)(-657,-200)
\put(-502,-131){\makebox(0,0)[lb]{\smash{{\SetFigFontNFSS{10}{12.0}{\familydefault}{\mddefault}{\updefault}$\mathcal B_1$}}}}
\put(830,1234){\makebox(0,0)[lb]{\smash{{\SetFigFontNFSS{10}{12.0}{\familydefault}{\mddefault}{\updefault}$\mathcal A_1$}}}}
\put(2650,-124){\makebox(0,0)[lb]{\smash{{\SetFigFontNFSS{10}{12.0}{\familydefault}{\mddefault}{\updefault}$\mathcal B_2$}}}}
\put(3100,1218){\makebox(0,0)[lb]{\smash{{\SetFigFontNFSS{10}{12.0}{\familydefault}{\mddefault}{\updefault}$\mathcal A_2$}}}}
\end{picture}%

\end{center}
\caption{Up to affine transformations, these are the only different pairs of convex cones with equal cross covariogram. }
\label{fig_settori_angolari}
\end{figure}
Clearly $\set{\F_1,-\G_1}\neq\set{\F_2,-\G_2}$. 
Elementary calculations prove $\jcov{\F_1,\G_1}=\jcov{\F_2,\G_2}$. Indeed, 
if  $x=(x_1,x_2)$ in Cartesian coordinates, then
\begin{equation*}
\jcov{\F_1,\G_1}(x)=\jcov{\F_2,\G_2}(x)=
\begin{cases}
x_2^2/2& \text{if $x\in \ang{0}{{\pi}/{4}}$;}\\
(x_2^2-x_1^2+2x_1x_2)/4& \text{if $x\in\ang{{\pi}/{4}}{{\pi}/{2}}$;}\\
(x_1+x_2)^2/4& \text{if $x\in\ang{{\pi}/{2}}{{3\pi}/{4}}$;}\\
0& \text{if $x\notin\ang{0}{{3\pi}/{4}}$.}
\end{cases}
\end{equation*}
These properties are preserved by any non-singular affine transformation $\cT$, since 
\[
\jcov{\cT\F_1,\cT\G_1}(x)=
|\det\cT|\jcov{\F_1,\G_1}(\cT^{-1}x)=|\det\cT|\jcov{\F_2,\G_2}(\cT^{-1}x)
=\jcov{\cT\F_2,\cT\G_2}(x).
\]
\end{example}

\begin{lemma}\label{discontinuita_piano}
Let $A$ and $B$ be pointed closed  convex cones in $\Real^2$ with non-empty interior satisfying $A$, $-B\subset\ang{0}{\pi}$. The set $S^2(A,B)=\cl\{x\in\Real^2 : g_{A,B}\text{ is not $C^2$ at $x$}\}$ coincides with $\pa A\cup(-\pa B)$.
\end{lemma}
\begin{proof}

Let $W\subset\Real^2\setminus\big(\pa A\cup(-\pa B)\big)$ be an open connected set. When $x\in W$ the vertex $x$ of $B+x$ does not belong to $\pa A$,and the vertex $O$ of $A$ does not belong to $\pa B+x$. Thus the combinatorial structure of $\pa(A\cap(B+x))$ is constant in $W$. Since the vertices of $A\cap(B+x)$ are smooth functions of $x$, for each  $x\in W$, so is $g_{A,B}$. This proves $S^2(A,B)\subset\pa A\cup(-\pa B)$.

Let $B=\ang{\al}{\be}$, for suitable $\al, \be\in[\pi,2\pi]$ with $\al-\be\neq\pm\pi$. It is easy to prove that, when $x\notin \pa A$, we have
\begin{align*}
\frac{\pa^2 g_{A,B}}{\pa u(\al)\pa u(\be)}(x)=&
-|\sin (\al-\be)|\,\frac{\pa}{\pa u(\al)}\len\big(A\cap(x+\{\te=\be\})\big)\\
=&-|\sin (\al-\be)|\, 1_A(x),
\end{align*}
This formula proves $\pa A\subset S^2(A,B)$. A similar formula for the second order mixed derivative of $g_{A,B}$ in the directions of the edges of $A$ proves $-\pa B\subset S^2(A,B)$.
\end{proof}

\begin{proof}[Proof of Theorem~\ref{cov_cong_settori_angolari}]We have $\conv(A\cup(-B))=\conv(A'\cup(-B'))=\supp g_{A,B}$, by \eqref{support}. Choose polar coordinates $(\rho,\te)$ so that $\supp g_{A,B}\subset\ang{0}{\pi}$.
Lemma~\ref{discontinuita_piano} proves that $\pa A\cup(-\pa B)$ is determined by $g_{A,B}$.
If $\pa A\cup(-\pa B)$ consists of two rays, then both $A$ and $-B$ coincide with the convex cone bounded by those rays. Therefore $\{A,-B\}$ is determined by $g_{A,B}$ and   $\jcov{A,B}=\jcov{A',B'}$ implies \eqref{acaso0}.

Assume that $\pa A\cup (-\pa B)$ consists of three rays.
Let  $0\leq\te_1<\te_2<\te_3\leq\pi$ be the angles corresponding to these rays.
Clearly exactly one among the rays that bound $A$ and $-B$ coincides with $\pa A\cap (-\pa B)$.
Elementary calculations show that, as $\ee\to0^+$,
\begin{equation*}
\jcov{A,B}(u(\te_1+\ee))=\begin{cases}
\ee+o(\ee) &\text{if $\{\te=\te_1\}=\pa A\cap(-\pa B)$};\\
o(\ee) &\text{otherwise.}
\end{cases}
\end{equation*}
An analogous formula, with $\{\te=\te_3\}$ substituting $\{\te=\te_1\}$, holds for $\jcov{A,B}(u(\te_3-\ee))$.
From the asymptotic behaviour of $\jcov{A,B}(u(\te_1+\ee))$ and $\jcov{A,B}(u(\te_3-\ee))$ it is thus possible to understand which of the tree rays  $\{\te=\te_1\}$, $\{\te=\te_2\}$ and $\{\te=\te_3\}$ coincides with $\pa A\cap(-\pa B)$.
If, for instance, $\{\te=\te_2\}\subset\pa A\cap(-\pa B)$, then   we necessarily have either
$A= \ang{\te_1}{\te_2}$ and $-B= \ang{\te_2}{\te_3}$ or else $-B= \ang{\te_1}{\te_2}$ and $A= \ang{\te_2}{\te_3}$. Thus  $\{A,-B\}$ is determined. Similar arguments prove that $\{A,-B\}$ is determined when $\pa A\cap(-\pa B)$ coincides with $\{\te=\te_1\}$ or with $\{\te=\te_3\}$.
The equality $\jcov{A,B}=\jcov{A',B'}$ implies \eqref{acaso0}.

Assume  that $\pa A\cup (-\pa B)$ consists of four rays,
say $\{\te=\te_i\}$, $i=1,\dots,4$, with  $0\leq\te_1<\te_2<\te_3<\te_4\leq\pi$. Let $P_2$ be the parallelogram bounded by the rays  $\{\te=\te_1\}$ and  $\{\te=\te_4\}$ and by the lines  which are parallel to these rays and contain $u(\te_2)$ (see Fig.~\ref{ang_sec_trefigure}). 
Let $a_2\neq O$ be the vertex of $P_2$ in $\{\te=\te_4\}$, $b_2\neq O$ be the vertex of $P_2$ in $\{\te=\te_1\}$ and let $c_2=\{\te=\te_3\}\cap[u(\te_2),a_2]$ .
There are three possible cases.

\begin{figure}
\begin{center}
\begin{picture}(0,0)%
\includegraphics{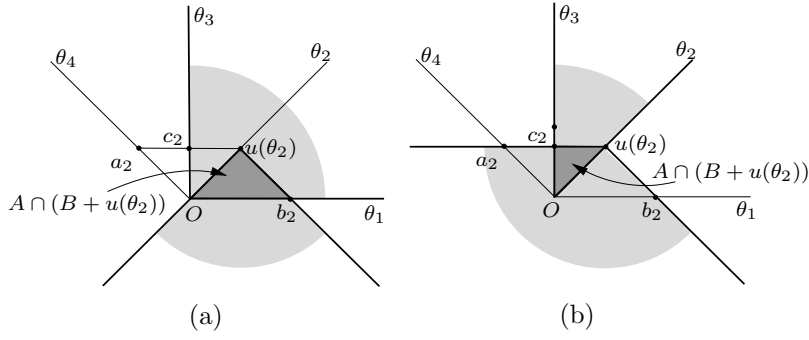}%
\end{picture}%
\setlength{\unitlength}{2901sp}%
\begingroup\makeatletter\ifx\SetFigFontNFSS\undefined%
\gdef\SetFigFontNFSS#1#2#3#4#5{%
  \reset@font\fontsize{#1}{#2pt}%
  \fontfamily{#3}\fontseries{#4}\fontshape{#5}%
  \selectfont}%
\fi\endgroup%
\begin{picture}(6304,2850)(1893,-324)
\put(4462,2022){\makebox(0,0)[lb]{\smash{{\SetFigFontNFSS{8}{9.6}{\familydefault}{\mddefault}{\updefault}{\color[rgb]{0,0,0}$\theta_2$}%
}}}}
\put(3477,2327){\makebox(0,0)[lb]{\smash{{\SetFigFontNFSS{8}{9.6}{\familydefault}{\mddefault}{\updefault}{\color[rgb]{0,0,0}$\theta_3$}%
}}}}
\put(3908,1202){\makebox(0,0)[lb]{\smash{{\SetFigFontNFSS{8}{9.6}{\familydefault}{\mddefault}{\updefault}{\color[rgb]{0,0,0}$u(\theta_2)$}%
}}}}
\put(4177,651){\makebox(0,0)[lb]{\smash{{\SetFigFontNFSS{8}{9.6}{\familydefault}{\mddefault}{\updefault}{\color[rgb]{0,0,0}$b_2$}%
}}}}
\put(2306,1969){\makebox(0,0)[lb]{\smash{{\SetFigFontNFSS{8}{9.6}{\familydefault}{\mddefault}{\updefault}{\color[rgb]{0,0,0}$\theta_4$}%
}}}}
\put(2772,1086){\makebox(0,0)[lb]{\smash{{\SetFigFontNFSS{8}{9.6}{\familydefault}{\mddefault}{\updefault}{\color[rgb]{0,0,0}$a_2$}%
}}}}
\put(4915,635){\makebox(0,0)[lb]{\smash{{\SetFigFontNFSS{8}{9.6}{\familydefault}{\mddefault}{\updefault}{\color[rgb]{0,0,0}$\theta_1$}%
}}}}
\put(3216,1304){\makebox(0,0)[lb]{\smash{{\SetFigFontNFSS{8}{9.6}{\familydefault}{\mddefault}{\updefault}{\color[rgb]{0,0,0}$c_2$}%
}}}}
\put(3433,-246){\makebox(0,0)[lb]{\smash{{\SetFigFontNFSS{10}{12.0}{\familydefault}{\mddefault}{\updefault}{\color[rgb]{0,0,0}(a)}%
}}}}
\put(6569,-241){\makebox(0,0)[lb]{\smash{{\SetFigFontNFSS{10}{12.0}{\familydefault}{\mddefault}{\updefault}{\color[rgb]{0,0,0}(b)}%
}}}}
\put(7550,2038){\makebox(0,0)[lb]{\smash{{\SetFigFontNFSS{8}{9.6}{\familydefault}{\mddefault}{\updefault}{\color[rgb]{0,0,0}$\theta_2$}%
}}}}
\put(6565,2343){\makebox(0,0)[lb]{\smash{{\SetFigFontNFSS{8}{9.6}{\familydefault}{\mddefault}{\updefault}{\color[rgb]{0,0,0}$\theta_3$}%
}}}}
\put(8063,672){\makebox(0,0)[lb]{\smash{{\SetFigFontNFSS{8}{9.6}{\familydefault}{\mddefault}{\updefault}{\color[rgb]{0,0,0}$\theta_1$}%
}}}}
\put(7265,667){\makebox(0,0)[lb]{\smash{{\SetFigFontNFSS{8}{9.6}{\familydefault}{\mddefault}{\updefault}{\color[rgb]{0,0,0}$b_2$}%
}}}}
\put(5394,1985){\makebox(0,0)[lb]{\smash{{\SetFigFontNFSS{8}{9.6}{\familydefault}{\mddefault}{\updefault}{\color[rgb]{0,0,0}$\theta_4$}%
}}}}
\put(5860,1102){\makebox(0,0)[lb]{\smash{{\SetFigFontNFSS{8}{9.6}{\familydefault}{\mddefault}{\updefault}{\color[rgb]{0,0,0}$a_2$}%
}}}}
\put(6297,1316){\makebox(0,0)[lb]{\smash{{\SetFigFontNFSS{8}{9.6}{\familydefault}{\mddefault}{\updefault}{\color[rgb]{0,0,0}$c_2$}%
}}}}
\put(7033,1233){\makebox(0,0)[lb]{\smash{{\SetFigFontNFSS{8}{9.6}{\familydefault}{\mddefault}{\updefault}{\color[rgb]{0,0,0}$u(\theta_2)$}%
}}}}
\put(3404,631){\makebox(0,0)[lb]{\smash{{\SetFigFontNFSS{8}{9.6}{\familydefault}{\mddefault}{\updefault}{\color[rgb]{0,0,0}$O$}%
}}}}
\put(6425,662){\makebox(0,0)[lb]{\smash{{\SetFigFontNFSS{8}{9.6}{\familydefault}{\mddefault}{\updefault}{\color[rgb]{0,0,0}$O$}%
}}}}
\put(1908,727){\makebox(0,0)[lb]{\smash{{\SetFigFontNFSS{8}{9.6}{\familydefault}{\mddefault}{\updefault}{\color[rgb]{0,0,0}$A\cap(B+u(\theta_2))$}%
}}}}
\put(7348,998){\makebox(0,0)[lb]{\smash{{\SetFigFontNFSS{8}{9.6}{\familydefault}{\mddefault}{\updefault}{\color[rgb]{0,0,0}$A\cap(B+u(\theta_2))$}%
}}}}
\end{picture}%
\end{center}
\caption{The set $A\cap(B+u(\theta_2))$ in Case~1 (left) and in Case~2 (right).}
\label{ang_sec_trefigure}
\end{figure}

\textit{Case 1. $\{A,-B\}=\{\ang{\te_1}{\te_3}, \ang{\te_2}{\te_4}\}$.}

It is easy to see that $A\cap(B+u(\te_2))$ is equal to the triangle $\conv\{O, u(\te_2), b_2\}$ when  $A= \ang{\te_1}{\te_3}$ and $-B= \ang{\te_2}{\te_4}$ (see Fig.~\ref{ang_sec_trefigure}(a)), and it is equal to the triangle $\conv\{O, u(\te_2), a_2\}$ when $A= \ang{\te_2}{\te_4}$ and $-B= \ang{\te_1}{\te_3}$. 
In each case we have
\begin{equation}\label{volumi_caso1}
\jcov{A,B}(u(\te_2))=\frac1{2}\area(P_2).
\end{equation}

\textit{Case 2.  $\{A,-B\}= \{\ang{\te_2}{\te_3}, \ang{\te_1}{\te_4}\}$.}

Suppose that $A= \ang{\te_2}{\te_3}$ and $-B= \ang{\te_1}{\te_4}$. The set $A\cap(B+u(\te_2))$ is strictly contained in the triangle $T=\conv(O, u(\te_2),a_2)$ (see Fig.~\ref{ang_sec_trefigure}(b)). 
Moreover the ratio between $\area(A\cap(B+u(\te_2)))$ and $\area(T)$ equals the ratio between $\|u(\te_2)-c_2\|$ and $\|u(\te_2)-a_2\|$.
Thus 
\begin{equation}\label{volumi_caso2a}
\jcov{A,B}(u(\te_2))
=\frac{\|u(\te_2)-c_2\|}{2\|u(\te_2)-a_2\|}\area(P_2)<\frac1{2}\area(P_2).
\end{equation}
The same formulas also hold when  $-B= \ang{\te_2}{\te_3}$ and $A= \ang{\te_1}{\te_4}$.

\textit{Case 3.  $\{A,-B\}=\{ \ang{\te_1}{\te_2}, \ang{\te_3}{\te_4} \}$.}

Arguments similar to those of Case~2 prove the following formula:
\begin{equation}\label{volumi_caso3a}
\jcov{A,B}(u(\te_2))
=\frac{\|c_2-a_2\|}{2\|u(\te_2)-a_2\|}\area(P_2)<\frac1{2}\area(P_2).
\end{equation}

The comparison of the values of $\jcov{A,B}(u(\te_2))$  in \eqref{volumi_caso1}, \eqref{volumi_caso2a} and \eqref{volumi_caso3a}  distinguishes Case~1 from the others.
Moreover, it distinguishes Case~2 from Case~3 except when $c_2$  divides the segment  $[u(\te_2),a_2]$  in two equal parts. Assume that this happens and let $\cT$ be a non-singular linear transformation  which maps the ray $\{\te=\te_i\}$ in $\{\te=(i-1)\pi/4\}$, for $i=1,2,4$. The assumption regarding $c_2$  easily implies $\cT\{\te=\te_3\}=\{\te=\pi/2\}$.

The same analysis, with the same $\te_i$, is also valid  for $A'$ and $B'$. 
If the same case applies to  $(A,B)$ and to $(A',B')$, then $\set{A,-B}=\set{A',-B'}$. 
The only possibility left is that there exists an affine transformation $\cT$ such that $\cT \{\te=\te_i\}=\{\te=(i-1)\pi/4\}$, for $i=1,\dots,4$, Case~2 applies to $(A,B)$ and Case~3 applies to $(A',B')$ (or vice versa). 
If this happens, then Alternative~\eqref{acaso12} in Theorem~\ref{cov_cong_settori_angolari} holds.
\end{proof}

\begin{remark}\label{vettore_in_comune} Observe that $\inte\F_1\cap\inte(-\G_1)\neq\emptyset$ and $\inte\F_2\cap\inte(-\G_2)=\emptyset$. Thus, if 
$\inte A\cap\inte(-B)$ and  $\inte A'\cap\inte (-B')$ are both empty or  both non-empty, then Theorem~\ref{cov_cong_settori_angolari} implies $\set{A,-B}=\set{A',-B'}$. 
\end{remark}
\begin{remark}\label{order_preserving} In the previous proof it is clear that the linear map $\cT$ in Theorem~\ref{cov_cong_settori_angolari} preserves the order of the rays, that is, the ray $\cT^{-1} \{\te=i\pi/4\}$ follows in counterclockwise order  the ray $\cT^{-1} \{\te=(i-1)\pi/4\}$, for $i=1,2,3$. 
\end{remark}

\section{Cross covariogram and synisothesis}\label{sec_crosscov_synisothesis} 
The next example is due to this author and R.~J.~Gardner.
\begin{example}\label{parall}
Let $\al,\be,\ga$ and $\de$ be positive real numbers, $I_1=[(-1,0),(1,0)]$, $I_2=1/\sqrt{2}\ [(-1,-1),(1,1)]$, $I_3=[(0,-1),(0,1)]$ and $I_4=1/\sqrt{2}\ [(1,-1),(-1,1)]$. We define four parallelograms as follows:
\begin{equation*}
\cK_1=\al I_1+\be I_2;\quad \cL_1=\ga I_3+\de I_4+y;\quad
\cK_2=\al I_1+\de I_4\quad \text{and}\quad \cL_2=\be I_2+\ga I_3+y.
\end{equation*}
\begin{figure}
\begin{center}
\begin{picture}(0,0)%
\includegraphics{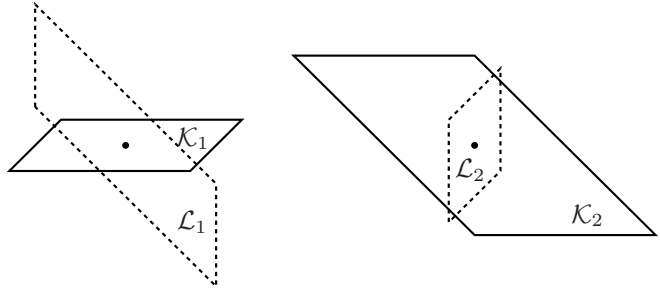}%
\end{picture}%
\setlength{\unitlength}{3522sp}%
\begingroup\makeatletter\ifx\SetFigFontNFSS\undefined%
\gdef\SetFigFontNFSS#1#2#3#4#5{%
  \reset@font\fontsize{#1}{#2pt}%
  \fontfamily{#3}\fontseries{#4}\fontshape{#5}%
  \selectfont}%
\fi\endgroup%
\begin{picture}(4544,2024)(-1191,-173)
\put(2746,299){\makebox(0,0)[lb]{\smash{{\SetFigFontNFSS{10}{12.0}{\familydefault}{\mddefault}{\updefault}$\mathcal K_2$}}}}
\put(1936,614){\makebox(0,0)[lb]{\smash{{\SetFigFontNFSS{10}{12.0}{\familydefault}{\mddefault}{\updefault}$\mathcal L_2$}}}}
\put(  1,839){\makebox(0,0)[lb]{\smash{{\SetFigFontNFSS{10}{12.0}{\familydefault}{\mddefault}{\updefault}$\mathcal K_1$}}}}
\put(  1,254){\makebox(0,0)[lb]{\smash{{\SetFigFontNFSS{10}{12.0}{\familydefault}{\mddefault}{\updefault}$\mathcal L_1$}}}}
\end{picture}%
\end{center}
\caption{Up to affine transformations, $(\cK_1,\cL_1)$ and $(\cK_2,\cL_2)$ are the only pairs of convex polygons with equal cross covariogram which are not synisothetic. }
\label{fig_parall}
\end{figure}
See Fig.~\ref{fig_parall}. The pairs $(\cK_1,-\cL_1)$ and $(\cK_2,-\cL_2)$ are not synisothetic (no vertex in the second pair has a support cone equal to the support cone of the top vertex of $\cL_1$ or to its reflection). Moreover $\jcov{\cK_1,\cL_1}=\jcov{\cK_2,\cL_2}$.
To prove  it,  let $X_i$ and $Y_i$ be independent random variables uniformly distributed over $\cK_i$ and $\cL_i$, for $i=1,2$, and let $Z_1,\dots,Z_4$ be independent random variables uniformly distributed over $I_1,\dots,I_4$, respectively.  Then $X_1=\al Z_1+\be Z_2$, $Y_1=\ga Z_3+\de Z_4+y$, $X_2=\al Z_1+\de Z_4$ and $Y_2=\be Z_2+\ga Z_3+y$.
Moreover we have
\begin{equation}\label{diff_rand_var}
X_1-Y_1=X_2-Y_2,
\end{equation}
because $Z_2=-Z_2$ and $Z_4=-Z_4$. The distribution of probability of $X_i-Y_i$ is $1_{\cK_i}\ast1_{-\cL_i}/(\area(\cK_i)\area(\cL_i))$. Observe that $\area(\cK_1)\area(\cL_1)=\al\be\ga\de/2=\area(\cK_2)\area(\cL_2)$. Therefore \eqref{diff_rand_var} and \eqref{convoluzione} imply
\begin{equation*}
g_{\cK_1,\cL_1}=1_{\cK_1}\ast1_{-\cL_1}=1_{\cK_2}\ast1_{-\cL_2}=g_{\cK_2,\cL_2}.
\end{equation*}
\end{example}

\begin{proposition}\label{det_curv}
Let $K$ and $L$ be convex polygons, $K'$ and $L'$ be planar closed convex sets with $g_{K,L}=g_{K',L'}$. Then $K'$ and $L'$ are polygons. 
Assume, moreover, that there is no affine transformation $\cT$ and no different indices $i,j\in\set{1,2}$ such that
$(\cT K,\cT L)$ and $(\cT K',\cT L')$ are trivial associates of $(\cK_i,\cL_i)$ and $(\cK_j,\cL_j)$, respectively. Then $(K,-L)$ and $(K',-L')$ are synisothetic.
\end{proposition}
This proof is divided in three steps and occupies all the rest of the section. We recall that $(K,-L)$ and $(K',-L')$ are synisothetic if and only if for each $u\in S^1$ both \eqref{sin_poligoni1} and \eqref{sin_poligoni2} hold. We will prove that \eqref{sin_poligoni1} holds for each $u$, while \eqref{sin_poligoni2} fails for some $u$ exactly when $(\cT K,\cT L)$ and $(\cT K',\cT L')$ are trivial associates of $(\cK_i,\cL_i)$ and $(\cK_j,\cL_j)$, respectively.

Observe that the set $K'-L'$ is a polygon, because $K'-L'=K-L$ by \eqref{support}. This may happen only if $K'$ and $L'$ are polygons.
The function $g_{K,L}$  determines its support $K-L$ and, for all $u\in S^1$, it determines also $\len(K_u)+\len(L_{-u})$, since \eqref{facce_corpodiff} implies
\begin{equation}\label{somma_lunghezze}
\len((K-L)_u)=\len(K_u)+\len(L_{-u}).
\end{equation}

\begin{claim}\label{step1} The function $g_{K,L}$ determines $\{\len(K_u), \len(L_{-u})\}$ for each $u\in S^1$. 
\end{claim} 

\begin{proof}If  $(K-L)_u$ is a vertex, then  both $K_u$ and $L_{-u}$ are vertices, by \eqref{somma_lunghezze}. In this case $\{\len(K_u), \len(L_{-u})\}=\{0\}$. Assume  that
\[
\text{$(K-L)_u$ is an edge.}
\]
In this case at least one among $K_u$ and $L_{-u}$ is an edge. Let  $x$ be the midpoint of $(K-L)_u$. For sufficiently small $\ee>0$,
let $x_{\ee}\in[O,x]$ be the point  at distance $\ee$ from $(K-L)_u$.
Then it is easy to see that, as $\ee\to0^+$, we have
\begin{equation}\label{order}
\jcov{K,L}(x_{\ee})=\min\{\len(K_u), \len(L_{-u})\}\ee+o(\ee).
\end{equation}
(This corresponds to translating $L$ so that the translated 
midpoint of $L_{-u}$ is close to the midpoint of $K_u$.)

If  $\jcov{K.L}(x_\ee)=o(\ee)$, then $\min\{\len(K_u), \len(L_{-u})\}=0$, that is, either $K_u$ is an edge and $L_{-u}$ is a vertex, or vice versa. If, for instance, $L_{-u}$ is a vertex, then the length of $K_u$ is determined, since it equals the length of $(K-L)_u$.

If  $\jcov{K.L}(x_\ee)\geq \al \,\ee$, for some constant $\al>0$, then both $K_u$ and $L_{-u}$ are edges. From \eqref{order} we obtain the minimum of the lengths of these edges.  Since the sum of these lengths is  determined by \eqref{somma_lunghezze}, the pair
$\{\len(K_u), \len(L_{-u})\}$ is known. \end{proof}

Let $z_0,\dots,z_n$, $w_0,\dots,w_t$, $z'_0,\dots,z'_{n'}$, $w'_0,\dots,w'_{t'}$ and $q_0,\dots,q_p$ denote respectively
the vertices of $K$, of $L$, of $K'$, of $L'$ and of $K-L=K'-L'$ in  counterclockwise order on the respective boundaries.
For each $j$, let $A_j=\cone(K,z_j)$, $B_j=\cone(L,w_j)$, $A'_j=\cone(K',z'_j)$ and $B'_j=\cone(L',w'_j)$.
Each vertex of $K-L$ is the difference of a vertex of $K$ and of one of $L$, and also of a vertex of $K'$ and of one of $L'$. Assume that $q_m=z_s-w_l=z'_h-w'_k$. In a neighbourhood of $z_s$, $K$ equals $A_s+z_s$, while in a neighbourhood of $w_l$, $L$ equals $B_l+w_l$. If $x$ belongs to a small neighbourhood of $q_m$, then we have
\[
K\cap(L+x)=(A_s+z_s)\cap(B_l+w_l+x)
\]
and this set is a translate of  $A_s\cap(B_l+x-q_m)$. The function $\jcov{A_s,B_l}$ is thus determined, in a neighbourhood of $O$, by $\jcov{K,L}$. 
Since $\jcov{A_s,B_l}$ is homogeneous of degree $2$ it is determined on its entire domain. Similar considerations apply to $(K',L')$ and imply $\jcov{A_s,B_l}=\jcov{A'_h,B'_k}$.

We say that the vertex $q_m=z_s-w_l=z'_h-w'_k$  is \emph{ambiguous}  if $\{A_s,-B_l\}\neq\{A'_h, -B'_k\}$. 

It is elementary to check that when $[z_s,z_{s+1}]$ and $[w_{l},w_{l+1}]$ are parallel we have 
\begin{equation}\label{lati_paralleli}
[q_m,q_{m+1}]= [z_s,z_{s+1}]-[w_{l},w_{l+1}],
\end{equation}
while when they are not parallel we have  either
\begin{align}
&\text{$[q_m,q_{m+1}]=[z_s,z_{s+1}]-w_l$ or }\label{edge_vertex}\\
&\text{$[q_m,q_{m+1}]=z_s-[w_{l},w_{l+1}]$.}\label{vertex_edge}
\end{align}

\begin{claim}\label{parallelismo} Let $q_m=z_s-w_l$ and assume that $q_m$ is ambiguous. Then $q_{m+1}$ is ambiguous too and \eqref{lati_paralleli} does not hold. Moreover, when  \eqref{edge_vertex} holds, $[z_{s-1},z_s]$ and $[z_{s+1},z_{s+2}]$ are parallel, while when \eqref{vertex_edge} holds, $[w_{l-1},w_{l}]$ and $[w_{l+1},w_{l+2}]$ are parallel.
\end{claim}

\begin{proof}Choose polar coordinates $(\rho,\theta)$ so that $[q_m,q_{m+1}]$ is parallel to $\{\theta=0\}$ and $A_s$ and $-B_l$  are contained in $\{0\leq \te\leq \pi\}$.  Since $\conv(A_s\cup(-B_l))=\conv(A'_h\cup(-B'_k))$, by \eqref{support}, also $A'_h$ and $-B'_k$ are contained in $\{0\leq \te\leq \pi\}$.
Theorem~\ref{cov_cong_settori_angolari} proves that there exist a linear transformation $\cT$ and $i$, $j\in\{1,2\}$, with $i\neq j$, such that Alternative~\eqref{acaso12} of Theorem~\ref{cov_cong_settori_angolari} occurs, with $A=A_s$, $B=B_l$, $A'=A'_h$ and $B'=B'_k$.
Since no edge of $\F_i$ is parallel to an edge of $-\G_i$, for each choice of $i$, \eqref{acaso1} implies that no edge of $K$ adjacent to $z_s$ is parallel to an edge of $L$ adjacent to $w_l$. This  rules out \eqref{lati_paralleli}. Analogous arguments rule out
$[q_m,q_{m+1}]= [z'_h,z'_{h+1}]-[w'_{k},w'_{k+1}]$. The latter implies that  $[q_m,q_{m+1}]$ equals either $[z'_h,z'_{h+1}]-w'_k$ or $z'_h-[w'_{k},w'_{k+1}]$.

Assume  \eqref{edge_vertex} 
and $[q_m,q_{m+1}]=[z'_h,z'_{h+1}]-w'_k$. This clearly implies
\begin{align}
A_s&= \ang{0}{\al},& &\phantom{ and }& A'_h&= \ang{0}{\al'},\label{apriori1}\\
A_{s+1}&= \ang{\be}{\pi}& &\text{ and }& A'_{h+1}&= \ang{\be'}{\pi},\label{apriori2}
\end{align}
for suitable $\al$, $\al'$, $\be$ and $\be'$ in $(0,\pi)$. It also implies $q_{m+1}=z_{s+1}-w_l=z'_{h+1}-w'_k$.
For $t=1,2,3,4$, let $\te_t\in[0,2\pi)$ be such that $\{\te=\te_t\}=\cT^{-1}\{\te=(t-1)\pi/4\}$. The cones $A_s$, $-B_l$, $A'_h$ and $-B'_k$ are bounded by the rays  $\{\te=\te_t\}$, by \eqref{acaso1}.  
By Remark~\ref{order_preserving}, and since these cones are contained in $\{0\leq\te\leq\pi\}$, we may assume
\begin{equation}\label{cond_ang1} 
0\leq\te_1<\te_2<\te_3<\te_4\leq\pi\quad\text{and}\quad\te_4\neq\te_1+\pi. 
\end{equation}
Moreover, the identities in \eqref{apriori1} imply that one of the $\te_t$, (necessarily $\te_1$, by \eqref{cond_ang1}) equals $0$. 
Assume $i=1$ and $j=2$ in \eqref{acaso1}. The condition \eqref{acaso1}, when expressed in terms of the $\te_t$,  becomes 
\begin{multline*}\label{aacaso1}
\{A_s,-B_l\}=\{ \ang{0}{\te_4}, \ang{\te_2}{\te_3}\}\text{ and}\\
\{A'_h,-B'_k\}=\{ \ang{0}{\te_2},\ang{\te_3}{\te_4}\}.
\end{multline*} Since \eqref{apriori1} implies $A_s\neq\ang{\te_2}{\te_3}$ and  $A_h\neq\ang{\te_3}{\te_4}$, we have
\begin{equation}\label{caso1}
A_s= \ang{0}{\te_4},\quad -B_l= \ang{\te_2}{\te_3}\quad\text{and}\quad -B'_k= \ang{\te_3}{\te_4}.
\end{equation}
Similar arguments imply that if $i=2$ and $j=1$ in \eqref{acaso1}, then we have
\begin{equation}\label{caso2}
A_s= \ang{0}{\te_2},\quad -B_l= \ang{\te_3}{\te_4}\quad\text{and}\quad -B'_k= \ang{\te_2}{\te_3}.
\end{equation}
Summarising, either \eqref{caso1} or \eqref{caso2} holds.

Let us prove that $q_{m+1}$ is ambiguous, that is, since $q_{m+1}=z_{s+1}-w_l=z'_{h+1}-w'_k$, let us   prove 
\begin{equation}\label{ambig_mpu}
\{A_{s+1},-B_l\}\neq \{A'_{h+1},-B'_k\}.
\end{equation}
The cone $-B'_k$ does not belong to the set in the left-hand side of \eqref{ambig_mpu}. Indeed, we have $\te_3<\te_4<\pi$, by \eqref{cond_ang1} and the equality $\te_1=0$. Thus \eqref{apriori2}, \eqref{caso1} and \eqref{caso2} imply $-B'_k\neq-B_l$ and $-B'_k\neq A_{s+1}$. 

Since $q_{m+1}$ is ambiguous, there exist a linear transformation $\calA$ and $i$, $j\in\{1,2\}$, with $i\neq j$, such that Alternative~\eqref{acaso12} of Theorem~\ref{cov_cong_settori_angolari} occurs, with $A=A_{s+1}$, $B=B_l$, $A'=A'_{h+1}$, $B'=B'_k$ and $\cT=\calA$. For $t=1,2,3,4$, let $\te'_t\in[0,2\pi)$ be such that $\{\te=\te'_t\}=\calA^{-1}\{\te=(t-1)\pi/4\}$. The $\te'_t$ satisfy a condition analogous to \eqref{cond_ang1} and, moreover, \eqref{apriori2} implies $\te'_4=\pi$.
It can be proved, by arguing as above, that  one of the  following possibilities occurs:
\begin{align}
A_{s+1}&=\ang{\te^\prime_1}{\pi},&  -B_l=&\ang{\te^\prime_2}{\te^\prime_3}\quad\text{and}&     -B'_k=&\ang{\te^\prime_1}{\te^\prime_2};\label{caso1p}\\
A_{s+1}&=\ang{\te^\prime_3}{\pi},&  -B_l=&\ang{\te^\prime_1}{\te^\prime_2}\quad\text{and}&     -B'_k=&\ang{\te^\prime_2}{\te^\prime_3}.\label{caso2p}
\end{align}
Observe that \eqref{caso1} and \eqref{caso1p} do not hold together, because \eqref{caso1} and \eqref{caso1p} imply $\te_3'=\te_3=\te_1'$, which contradicts $\te_3'>\te_1'$. Similar arguments prove that \eqref{caso2} and \eqref{caso2p} do not hold together. 
Assume \eqref{caso1} and \eqref{caso2p}.
In this case we have $\te_4=\te_3'$. Thus the identities in \eqref{caso2} and \eqref{caso2p} regarding $A_s$ and $A_{s+1}$ become 
\[
A_s=\ang{0}{\te_4}\quad\text{and}\quad A_{s+1}=\ang{\te_4}{\pi}.
\]
These conditions clearly imply that $[z_{s-1},z_s]$ and $[z_{s+1},z_{s+2}]$ are parallel to $\{\te=\te_4\}$. When \eqref{caso2} and \eqref{caso1p} hold we have $\te_2=\te_1'$. Thus we have $A_s=\ang{0}{\te_2}$ and $A_{s+1}=\ang{\te_2}{\pi}$, 
which again imply  $[z_{s-1},z_s]$ parallel to $[z_{s+1},z_{s+2}]$. 
This concludes the proof if \eqref{edge_vertex} holds and $[q_m,q_{m+1}]=[z'_h,z'_{h+1}]-w'_k$. If  \eqref{edge_vertex} holds and  $[q_m,q_{m+1}]=z'_h-[w'_k,w'_{k+1}]$, then the claim can be proved as before, by substituting in the proof $A'_h$, $A'_{h+1}$ and $-B'_k$ respectively with $-B'_k$, $-B'_{k+1}$ and $A'_h$. Similar arguments prove the claim if \eqref{edge_vertex} is substituted by \eqref{vertex_edge}. 
\end{proof}

\begin{claim}\label{step2} The function $g_{K,L}$ determines $\set{\cone(K,K_u),\cone(-L,(-L)_u)}$ for each $u\in S^1$.
\end{claim}

\begin{proof}Assume that no vertex of $K-L$ is ambiguous.
Let $u\in S^1$. Claim~\ref{step1} implies that $g_{K,L}$ distinguishes whether both $K_u$ and $L_{-u}$ are vertices,  or both $K_u$ and $L_{-u}$ are edges or one is a vertex and the other one is an edge.
If both $K_u$ and $L_{-u}$ are vertices, then the claim follows from the assumption that $(K-L)_u$ is not ambiguous.
If $K_u$ and $L_{-u}$ are edges, then $\set{\cone(K,K_u), \cone(-L,(-L)_u)}=\{H\}$, where $H=\{x\in\Real^2:x\cdot u\leq0\}$. 
If $K_u$ is an edge  and $L_{-u}$ is a vertex (or vice versa), then consider the set
\[
\bigcup_{w\in S^1 : (K-L)_w\text{ is a vertex}}\{\cone(K,K_w),\cone(-L,(-L)_w)\}.
\]
This set is determined by $g_{K,L}$, since no vertex of $K-L$ is ambiguous. The convexity of $K$ and of $L$ implies that only one among these cones, say $A$, has the property that $A\setminus\{O\}\subset \inte H$.  Then 
\[
 \set{\cone(K,K_u),\cone(-L,(-L)_u)}=\set{H, A}.
\]

Assume  that some vertex of $K-L$ is ambiguous. By Claim~\ref{parallelismo}, all vertices are ambiguous. Let us use the notations introduced before Claim~\ref{parallelismo}.
Let  $s$ be any index in $\set{1,\dots,n}$ and choose $l\in\set{1,\dots,t}$ in such a way that $[z_s,z_{s+1}]-w_l$ is an edge of $K-L$. 
This is possible because if $u$ denotes  the outer normal to $K$ at $[z_s,z_{s+1}]$, then  $L_{-u}$ is a vertex, by Claim~\ref{parallelismo}. 
By the same claim, $[z_{s-1},z_s]$ and $[z_{s+1},z_{s+2}]$  are parallel.
A similar argument proves that, given any $l\in\set{1,\dots,t}$, $[w_{l-1},w_{l}]$ and $[w_{l+1},w_{l+2}]$  are parallel. 
This is possible only if both $K$ and $L$ are parallelograms.
Similar arguments prove that $K'$ and $L'$ are parallelograms too.

Consider a given vertex $q_m=z_s-w_l=z'_h-w'_k$ of $K-L$.
Since $q_m$ is ambiguous there exists a linear transformation $\cT$ such that Alternative~\eqref{acaso12} of Theorem~\ref{cov_cong_settori_angolari} holds, with $A_s=A$, $B_l=B$, $A'_h=A'$ and $B'_k=B'$.
Assume, for instance, 
\[
\cT A_s=\F_1,\quad\cT B_l=\G_1,\quad\cT A'_h=\F_2\quad\text{and}\quad\cT B_l=\G_2.
\]
Since $\cT K$ is a parallelogram, the set of the directions of its edges  coincides with the set of the directions of the edges of  the support cone $\cT A_s$ in one  vertex of $\cT K$. Thus the edges of $\cT K$ are parallel to those of $\cK_2$.
Similarly, the edges of $\cT L$ (of  $\cT K'$ and $\cT L'$) 
are parallel to those of $\cL_2$ (of $\cK_1$ and $\cL_1$, respectively). Let $x_1$ and $x_2$ be the centres of $\cT K$ and of $\cT K'$, respectively, and let $y$ be the center of $\cT L-x_1$ and $\cT L'-x_2$ (these sets have equal center because $K-L=K'-L'$, by \eqref{support}). Choose the parameters defining $\cK_2$ and $\cL_2$ so that $\cT K-x_1=\cK_2$ and $\cT L-x_1=\cL_2+y$.
The edges of $\cT K$ and $\cT K'$ parallel to $\{\te=0\}$ have equal length, by Claim~\ref{step1} and because $\cT L$ and $\cT L'$ have no edges parallel to $\{\te=0\}$. Also the edges of $\cT K$ and $\cT L'$ parallel to $\{\te=3\pi/4\}$ have equal length, and the same property holds for the edges of $\cT K'$ and $\cT L$ parallel to $\{\te=\pi/4\}$, and for those of $\cT L$ and $\cT L'$ parallel to $\{\te=\pi/2\}$. Therefore $\cT K'-x_2=\cK_1$ and $\cT L'-x_2=\cL_1+y$. 
This contradicts the assumptions of Proposition~\ref{det_curv} and proves the claim.
Claim~\ref{step1} and Claim~\ref{step2} imply that $(K,-L)$ and $(K',-L')$ are synisothetic and conclude the proof of Proposition~\ref{det_curv}.
\end{proof}


\section{A crucial lemma}\label{sec_crucial_lemma}
This section is dedicated to the proof of Lemma~\ref{archi_simmetrici} and to some results needed in its proof. 
The first one, Lemma~\ref{coni_finale}, is, in our opinion, of interest by itself. It is contained in the unpublished note~\cite{ML}, where it is proved with geometrical arguments. 
Here we present a different, shorter proof which is based on the Theorem of supports for convolutions~\cite[Th.~4.3.3]{H}. 
\begin{lemma}\label{coni_finale}
Let $A$, $B$, $C$ and $D$ be convex cones in  $\Real^n$, $n\geq2$, with apex the origin $O$. Assume that each of them either coincides with $\set{O}$ or has non-empty interior and, moreover,
$A\cup B\subset \set{(x_1,x_2,\dots,x_n)\in\Real^n:x_n\geq0}$,
$A\cap\set{x_n=0}=B\cap\set{x_n=0}=\set{O}$,
$C\cup D\subset \set{x_n\leq0}$ and
$\conv (C\cup D)$ is pointed.
If
\begin{equation}\label{due_somme}
g_{A,C}+g_{B,D}=g_{A,D}+g_{B,C}
\end{equation}
then either $A=B$ or  $C=D$. The same conclusion holds if the  hypothesis ``$\conv (C\cup D)$  is pointed'' is substituted by ``either $A\subset B$ or $B\subset A$, and, moreover, either $C\subset D$ or $D\subset C$''.
\end{lemma}
\begin{proof}
For $r>0$, let $A_r=A\cap B(O,r)$, $B_r=B\cap B(O,r)$, $C_r=C\cap B(O,r)$ and $D_r=D\cap B(O,r)$.
We prove that there exists $s>1$ such that, for each $x\in B(O,1)$, we have
\begin{equation}\label{A_s}\begin{aligned}
g_{A,C}(x)&=g_{A_s,C_s}(x),& g_{B,D}(x)&=g_{B_s,D_s}(x),\\ g_{B,C}(x)&=g_{B_s,C_s}(x)\text{ and}& g_{A,D}(x)&=g_{A_s,D_s}(x).
\end{aligned}\end{equation}
Since $A\cap\{x_n=0\}=B\cap\{x_n=0\}=\{O\}$, there exists $s>1$ such that $(A\cup B)\cap\set{x_n\leq1}\subset B(O,s-1)$.
Let $x\in B(O,1)$. Since $C+x\subset\set{x_n\leq1}$, we have $A\cap(C+x)\subset A\cap\set{x_n\leq1}\subset B(O,s-1)$. Moreover, by the triangle inequality, we have $(C+x)\cap B(O,s-1)\subset(C_s+x)\cap B(O,s-1)$. Therefore $A\cap(C+x)=A_s\cap(C_s+x)$ and $g_{A,C}(x)=g_{A_s,C_s}(x)$.
Similar arguments prove the other identities in \eqref{A_s}.

All the functions which appear in \eqref{due_somme} are homogeneous of degree $n$ and  \eqref{due_somme} holds true if and only if it holds true in $B(O,1)$, that is, if and only if
$g_{A_s,C_s}(x)+g_{B_s,D_s}(x)=g_{A_s,D_s}(x)+g_{B_s,C_s}(x)$ for each  $x\in B(O,1)$.
By  \eqref{convoluzione} this condition  is equivalent to
\begin{equation}\label{diff_conv}
(1_{A_s}-1_{B_s})*(1_{-C_s}-1_{-D_s})(x)=0\quad
\text{for each  $x\in B(O,1)$.}
\end{equation}

Let us conclude the proof under the assumption $\conv (C\cup D)$  pointed.
Let $S= \supp (1_{A_s}-1_{B_s})*(1_{-C_s}-1_{-D_s})$. By \eqref{diff_conv}, we have $S\cap B(O,1)=\emptyset$. The set  $S$ is clearly contained in $\conv (A\cup B\cup (-C)\cup(-D))$ and the assumptions of the lemma imply that this union is pointed.  Therefore the identity $S\cap B(O,1)=\emptyset$ implies that there exists $\ee>0$ such that
\begin{equation*}\label{cono_non_interseca}
(\conv S)\cap B(O,\ee)=\emptyset.
\end{equation*}
We may apply the Theorem of supports for convolutions~\cite[Th.~4.3.3]{H}, since the  involved functions have compact supports. This theorem implies
\begin{equation*}\label{supp_conv}
\conv\ S=\conv \supp (1_{A_s}-1_{B_s})+\conv \supp (1_{-C_s}-1_{-D_s}).
\end{equation*}
Therefore either we have $\conv \supp (1_{A_s}-1_{B_s})\cap B(O,\ee/2)=\emptyset$ or  we have $\conv \supp (1_{-C_s}-1_{-D_s})\cap B(O,\ee/2)=\emptyset$. In the first case we have $A_s\cap B(O,\ee/2)=B_s\cap B(O,\ee/2)$, which is equivalent to $A=B$. In the second case, by similar arguments, we have $C=D$.

Drop the assumption $\conv (C\cup D)$  pointed. When $A\subset B$ and $C\subset D$ then the functions which are convolved in \eqref{diff_conv} are constant in the interior of their supports, and their convolution vanish in $B(O,1)$ if and only if one of them vanish. The other cases are treated similarly.
\end{proof}
Let us introduce the second counterexample for the cross covariogram problem for convex polygons.
\begin{example}\label{parall_due}Let $\al,\be,\ga$ and $\de$ be positive real numbers, $m\in \Real$, $y\in\Real^2$, $I_1$ and $I_3$ be as in Example~\ref{parall} and  $I^{(m)}=(1/\sqrt{1+m^2})\,[(-m,-1),(m,1)]$. Assume either  $m=0$, $\al\neq \ga$ and $\be\neq\de$ or else $m\neq 0$ and $\al\neq \ga$. We define four parallelograms as follows:
\begin{equation*}
\cK_3=\al I_1+\be I_3;\quad
\cL_3=\ga I_1+\de I^{(m)}+y;\quad
\cK_4=\ga I_1+\be I_3;\quad\text{and}\quad
\cL_4=\al I_1+\de I^{(m)}+y.
\end{equation*}
See~Fig.~\ref{fig_parall_due}. We have $\jcov{\cK_3,\cL_3}=\jcov{\cK_4,\cL_4}$ (it can be proved by arguing as in Example~\ref{parall}), and the pairs $(\cK_3,-\cL_3)$ and $(\cK_4,-\cL_4)$ are clearly synisothetic. 
However, $(\cK_3,\cL_3)$ and $(\cK_4,\cL_4)$ are not trivial associates. 

\begin{figure}
\begin{center}
\begin{picture}(0,0)%
\includegraphics{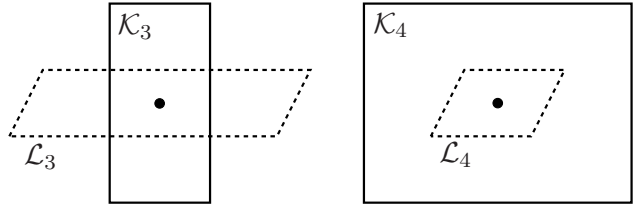}%
\end{picture}%
\setlength{\unitlength}{3646sp}%
\begingroup\makeatletter\ifx\SetFigFontNFSS\undefined%
\gdef\SetFigFontNFSS#1#2#3#4#5{%
  \reset@font\fontsize{#1}{#2pt}%
  \fontfamily{#3}\fontseries{#4}\fontshape{#5}%
  \selectfont}%
\fi\endgroup%
\begin{picture}(4229,1394)(-831,817)
\put(-79,1980){\makebox(0,0)[lb]{\smash{{\SetFigFontNFSS{11}{13.2}{\familydefault}{\mddefault}{\updefault}$\mathcal K_3$}}}}
\put(-710,1101){\makebox(0,0)[lb]{\smash{{\SetFigFontNFSS{11}{13.2}{\familydefault}{\mddefault}{\updefault}$\mathcal L_3$}}}}
\put(1641,1992){\makebox(0,0)[lb]{\smash{{\SetFigFontNFSS{11}{13.2}{\familydefault}{\mddefault}{\updefault}$\mathcal K_4$}}}}
\put(2083,1113){\makebox(0,0)[lb]{\smash{{\SetFigFontNFSS{11}{13.2}{\familydefault}{\mddefault}{\updefault}$\mathcal L_4$}}}}
\end{picture}%
\end{center}
\caption{Up to affine transformations, $(\cK_3,\cL_3)$ and $(\cK_4,\cL_4)$ are the only pairs of convex polygons with equal cross covariogram which are synisothetic and are not trivial associates.}
\label{fig_parall_due}
\end{figure}
\end{example}

\begin{lemma}\label{coincidono_intorno_lato}
Let $K$, $L$, $K'$ and $L'$ be convex polygons satisfying the assumptions of Theorem~\ref{cov_congiunto_poligoni} with $(K,-L)$ and $(K',-L')$ synisothetic.
Assume that, for a given $u\in S^1$, $K_u$ and $K'_u$ are edges and have different lengths. Let $I^1_K$ (and $I^2_K$) be the edge of $K$ adjacent to $K_u$ that, in counterclockwise order on $\pa K$, precedes (and follows, respectively) $K_u$. 
For $i=1,2$, let $I^i_{-L}$, $I^i_{K'}$ and $I^i_{-L'}$  be respectively edges of $-L$, of $K'$ and of $-L'$ defined in analogy to $I^i_K$. Then $(-L)_u$ and $(-L')_u$ are edges and 
\begin{enumerate}
\item\label{coincidono1}  either $I^i_K$ is parallel to $I^i_{-L'}$ and $I^i_{K'}$ is parallel to $I^i_{-L}$, for $i=1,2$, 
\item\label{coincidono2}  or $I^1_K$, $I^2_K$, $I^1_{K'}$ and $I^1_{K'}$ are parallel  and $I^1_{-L}$, $I^2_{-L}$, $I^1_{-L'}$ and $I^2_{-L'}$ are parallel.
\end{enumerate} 
\end{lemma}

\begin{remark}When $(K,L)=(\cK_3,\cL_3)$, $(K',L')=(\cK_4,\cL_4)$ and $u=(0,1)$, then Alternative~\eqref{coincidono2} of Lemma~\ref{coincidono_intorno_lato} occurs; see Fig.~\ref{fig_parall_due}.
\end{remark}

\begin{proof}
Since $\len(K_u)\neq\len(K'_u)$, the synisothesis of $(K,-L)$ and $(K',-L')$ implies  \eqref{alt_due}. Therefore $(-L)_u$ and $(-L')_u$ are edges, $\len (K'_u)=\len ((-L)_u)$ and $\len((-L')_u)=\len (K_u)$. 
Let  $K_u=[z_0, z_1]$ and $(-L)_u=[w_0, w_1]$, where $z_1$ follows $z_0$ and $w_1$ follows $w_0$, in counterclockwise order on the respective boundaries. 
For $i=1,2$, let $v^i_K\in S^1$ be parallel to $I^i_K$ and oriented in such a way that $v^i_K\cdot u>0$.
Define $v^i_{-L}$, $v^i_{K'}$ and  $v^i_{-L'}$ similarly.
Assume, for instance, $\len (K_u)>\len(K'_u)$, that is, 
$\len (K_u)>\len((-L)_u)$ and $\len (K'_u)<\len((-L')_u)$. Let $q_0=z_0+w_0$, $q_1=z_0+w_1$, $q_2=z_1+w_0$ and $q_3=z_1+w_1$. 
The points  $q_0$, $q_1$, $q_2$ and $q_3$  belong to $[q_0,q_3]=K_u+(-L)_u$ Moreover we have $[q_0,q_3]=(K-L)_u$ (by \eqref{facce_corpodiff}) and  $q_0<q_1<q_2<q_3$ in counterclockwise order on $\pa(K-L)$.

Let $S^2(K,L)=\cl\{x\in \Real^2 : \text{$g_{K,L}$ is not $C^2$  at $x$}\}$. We analyse the shape of $S^2(K,L)\cap W$, where $W$ is a  neighbourhood of $[q_0,q_3]$. It is easy to prove that 
\[
S^2(K,L)=\Big(\cup_{\text{$z$ vertex of $K$}}(-\pa L+z)\Big)\cup\Big(\cup_{\text{$w$ vertex of $-L$}}(\pa K+w)\Big).
\] 
Schmitt~\cite{S} proves this formula when $K=L$ and the general case can be proved in the same way. We also recall Lemma~\ref{discontinuita_piano}, which proves the previous formula when $K$ and $L$ are planar convex cones. Thus, when $W$ is sufficiently small, we have 
\[
W\cap S^2(K,L)=W\cap\Big((-\pa L+z_0)\cup(-\pa L+z_1)\cup(\pa K+w_0)\cup(\pa K+w_1)\Big).
\]
This set is the union  of $[q_0,q_3]$ and, for each $i=0,\dots,3$, of two  line segments (possibly coincident) containing $q_i$.
If $U_i$ denotes the set of the directions of the line segments containing $q_i$, then
\[ 
U_0=\{v^1_K,v^1_{-L}\},\ U_1=\{v^1_K,v^2_{-L}\}, \ U_2=\{v^2_K,v^1_{-L}\}\text{ and } U_3=\{v^2_K,v^2_{-L}\}.
\]
The above analysis can be repeated for $g_{K',L'}$. However,  in this case $K'$ has the role of $-L$ and $-L'$ the role of $K$, because $\len (K'_u)<\len((-L)'_u)$. Therefore the identity $g_{K,L}=g_{K',L'}$ implies
\[
U_0=\{v^1_{K'},v^1_{-L'}\},\ U_1=\{v^2_{K'},v^1_{-L'}\}, \ U_2=\{v^1_{K'},v^2_{-L'}\}\text{ and } U_3=\{v^2_{K'},v^2_{-L'}\}.
\]

Observe that if one of the equalities $v^1_K=v^1_{-L'}$, $v^2_K=v^2_{-L'}$, $v^1_{-L}=v^1_{K'}$ and $v^2_{-L}=v^2_{K'}$  holds, then also the other three equalities hold. 
For instance, if $v^1_K=v^1_{-L'}$, then the identities involving $U_0$ imply $v^1_{-L}=v^1_{K'}$, those involving  $U_1$ imply $v^2_{-L}=v^2_{K'}$.  
Once these are established, the identities involving $U_2$ imply $v^2_K=v^2_{-L'}$. When one of these equalities holds, Alternative~\eqref{coincidono1} occurs.

Assume  $U_0=U_1=U_2=\{v,w\}$, for suitable $v,w$, with $v\neq w$. 
If $v^1_{K}=v^1_{-L'}$, then Alternative~\eqref{coincidono1} occurs, as proved above. 
If $v^1_{K}\neq v^1_{-L'}$ and, for instance, $v^1_{K}=w$ and $v^1_{-L'}=v$, then necessarily $v^1_{K}=v^2_{K}=w$, $v^1_{-L}=v^2_{-L}=v$, $v^1_{K'}=v^2_{K'}=w$ and $v^1_{-L'}=v^2_{-L'}=v$, that is, Alternative~\eqref{coincidono2} occurs. 
When $U_0=U_1=U_2$ consists of a single element, clearly $v^1_{K}=v^1_{-L'}$ and Alternative~\eqref{coincidono1} occurs.
If $U_0\neq U_1$, then $v^1_{K}$ and $v^1_{-L'}$ coincide, because they are the only element of $U_0\cap U_1$, and Alternative~\eqref{coincidono1} occurs. Finally,
if $U_0\neq U_2$, then $v^1_{K'}$ and $v^1_{-L}$ coincide, because they are the only element of $U_0\cap U_2$, and again Alternative~\eqref{coincidono1} occurs. 
\end{proof}

\begin{lemma}\label{archi_simmetrici}
Assume that $K$, $K'$, $L$ and $L'$ are convex polygons satisfying the assumptions of Theorem~\ref{cov_congiunto_poligoni}, that $(K,-L)$ and $(K',-L')$ are synisothetic and that $(K,L)$ and $(K',L')$ are not trivial associates. Assume also the following properties:
\begin{enumerate}
\item\label{ass1} there exists an arc $U\subset S^1$, which is not a point, such that \eqref{alt_uno} holds true for each $u\in U$, and $U$ is a maximal arc (with respect to inclusion) with this property;
\item\label{ass2}
there exist $u_0\in U$ such that  $K_{u_0}=K'_{u_0}$ and
$(-L)_{u_0}=(-L')_{u_0}$;
\item\label{ass3}  if $\Si$ denotes the maximal closed arc contained in $\partial K\cap\partial K'$ and containing $K_{u_0}$, and  $\Om$ denotes the maximal closed arc contained in $\partial (-L)\cap \partial(- L')$ and containing $(-L)_{u_0}$, then neither $\Si$ nor $\Om$ are points or line segments.
\end{enumerate}
Then $\Si$ is a translate of $\Om$. \end{lemma}

The proof of this lemma  is divided in five steps and occupies all the rest of the section. 
First observe that  $\Si\neq\pa K$ and $\Om\neq\pa(-L)$. 
Indeed, if, for instance, $\Si=\pa K$, then $K=K'$. Moreover $L=L'$, because  $K-L=K'-L'$ (by \eqref{support}) and  the Minkowski addition satisfies a cancellation law (see~\cite[p.~126]{Sc}). Thus $(K,L)$ and $(K',L')$ are trivial associates, a contradiction. 

\begin{claim}\label{claim1}
Let $a_1$, $a_2\in\pa K$, $b_1$, $b_2\in\pa (-L)$ and $u_1$, $u_2\in S^1$ satisfy 
$\Si=[a_1,a_2]_{\pa K}$, $\Om=[b_1,b_2]_{\pa(- L)}$ and $\cl U=[u_1,u_2]_{S^1}$; see Fig.~\ref{archi_setting}.
Then, for each $j=1,2$, the sets $K_{u_j}$ and $K'_{u_j}$ are edges which contain $a_j$ and have a common  endpoint $a_j'$ contained in $\relint \Si$. 
Similarly, the sets $(-L)_{u_j}$ and $(-L')_{u_j}$ are edges which contain $b_j$ and have a common endpoint $b_j'$ contained in $\relint \Om$.
In particular, $\Si$ and $\Om+a_j-b_j$ coincide in a neighbourhood of $a_j$.

Moreover, if, for some $j\in\{1,2\}$, either $\len(K_{u_j})\neq\len(K'_{u_j})$  and Alternative~\eqref{coincidono1} of Lemma~\ref{coincidono_intorno_lato} (with $u$ replaced by $u_j$) holds, or else $\len(K_{u_j})=\len(K'_{u_j})$, then
\begin{equation}\label{coni_giusti}
N(K,a_j)= N(-L',b_j)\quad\text{and}\quad N(K',a_j)= N(-L,b_j).
\end{equation}
\end{claim}

\begin{figure}
\begin{center}
\begin{picture}(0,0)%
\includegraphics{archi_setting.pstex}%
\end{picture}%
\setlength{\unitlength}{3108sp}%
\begingroup\makeatletter\ifx\SetFigFontNFSS\undefined%
\gdef\SetFigFontNFSS#1#2#3#4#5{%
  \reset@font\fontsize{#1}{#2pt}%
  \fontfamily{#3}\fontseries{#4}\fontshape{#5}%
  \selectfont}%
\fi\endgroup%
\begin{picture}(6384,1642)(-2180,1170)
\put(-1053,1680){\makebox(0,0)[lb]{\smash{{\SetFigFontNFSS{9}{10.8}{\familydefault}{\mddefault}{\updefault}$u_2$}}}}
\put(-2165,1713){\makebox(0,0)[lb]{\smash{{\SetFigFontNFSS{9}{10.8}{\familydefault}{\mddefault}{\updefault}$u_1$}}}}
\put(-1696,1295){\makebox(0,0)[lb]{\smash{{\SetFigFontNFSS{9}{10.8}{\familydefault}{\mddefault}{\updefault}$U$}}}}
\put(1266,1988){\makebox(0,0)[lb]{\smash{{\SetFigFontNFSS{9}{10.8}{\familydefault}{\mddefault}{\updefault}$a_2$}}}}
\put(-494,1882){\makebox(0,0)[lb]{\smash{{\SetFigFontNFSS{9}{10.8}{\familydefault}{\mddefault}{\updefault}$a_1$}}}}
\put(1274,2471){\makebox(0,0)[lb]{\smash{{\SetFigFontNFSS{9}{10.8}{\familydefault}{\mddefault}{\updefault}$\pa K'$}}}}
\put(870,2466){\makebox(0,0)[lb]{\smash{{\SetFigFontNFSS{9}{10.8}{\familydefault}{\mddefault}{\updefault}$\pa K$}}}}
\put(2477,1248){\makebox(0,0)[lb]{\smash{{\SetFigFontNFSS{9}{10.8}{\familydefault}{\mddefault}{\updefault}$\Omega$}}}}
\put(1877,1894){\makebox(0,0)[lb]{\smash{{\SetFigFontNFSS{9}{10.8}{\familydefault}{\mddefault}{\updefault}$b_1$}}}}
\put(3751,2193){\makebox(0,0)[lb]{\smash{{\SetFigFontNFSS{9}{10.8}{\familydefault}{\mddefault}{\updefault}$b_2$}}}}
\put(3260,2635){\makebox(0,0)[lb]{\smash{{\SetFigFontNFSS{9}{10.8}{\familydefault}{\mddefault}{\updefault}$\pa(-L')$}}}}
\put(3974,2641){\makebox(0,0)[lb]{\smash{{\SetFigFontNFSS{9}{10.8}{\familydefault}{\mddefault}{\updefault}$\pa(-L)$}}}}
\put(1171,1739){\makebox(0,0)[lb]{\smash{{\SetFigFontNFSS{9}{10.8}{\familydefault}{\mddefault}{\updefault}$a_2'$}}}}
\put(-44,1244){\makebox(0,0)[lb]{\smash{{\SetFigFontNFSS{9}{10.8}{\familydefault}{\mddefault}{\updefault}$\Sigma$}}}}
\put(-404,1424){\makebox(0,0)[lb]{\smash{{\SetFigFontNFSS{9}{10.8}{\familydefault}{\mddefault}{\updefault}$a_1'$}}}}
\put(2026,1379){\makebox(0,0)[lb]{\smash{{\SetFigFontNFSS{9}{10.8}{\familydefault}{\mddefault}{\updefault}$b_1'$}}}}
\put(3871,1694){\makebox(0,0)[lb]{\smash{{\SetFigFontNFSS{9}{10.8}{\familydefault}{\mddefault}{\updefault}$b_2'$}}}}
\end{picture}%
\end{center}
\caption{The sets $U$, $\Si$ and $\Om$.}
\label{archi_setting}
\end{figure}

\begin{proof}
Let $u^0_1$ be the upper endpoint of $S^1\cap N(K,a_1)$ and of $S^1\cap N(K',a_1)$. These endpoints coincide because $[a_1,a_2]_{\pa K}\subset\pa K\cap \pa K'$. Define $u^0_2$ as $u^0_1$, with lower replacing upper and $a_2$ replacing $a_1$. Since $K_{u_0}=K'_{u_0}\subset\Si$ and $K_{u_0}$ is isothetic to $K'_{u_0}$, by assumption, we have $u_0\in[u^0_1,u^0_2]_{S^1}$. Let us prove
\begin{equation}\label{caratt_U}
 [u^0_1,u^0_2]_{S^1}=\cl U.
\end{equation}
If $u\in (u^0_1,u^0_2)_{S^1}$ then $K_u$ and $K'_u$ are isothetic, because $K_u$ and $K'_u$ are contained in $\relint\Si$ and therefore $\pa K$ and $\pa K'$ coincide in a neighbourhood of $K_u=K'_u$.  The synisothesis of $(K,-L)$ and $(K',-L')$ implies that 
also $(-L)_u$ and $(-L')_u$ are isothetic. Thus \eqref{alt_uno} holds for each $u\in(u^0_1,u^0_2)_{S^1}$. Since $[u^0_1,u^0_2]_{S^1}$ intersects $U$ (both arcs contains $u_0$) and $U$ is maximal, we have  $[u^0_1,u^0_2]_{S^1}\subset\cl U$.

In order to conclude the proof of \eqref{caratt_U} it suffices to show that  in any neighbourhood of $u_j^0$ there are directions $u$ for which  \eqref{alt_uno} does not hold, for $j=1,2$. Since $\pa K$ and $\pa K'$ bifurcate at $a_j$, $a_j$ is a vertex of $K$ or of $K'$. If it is both a vertex of $K$ and of $K'$ then the vertex $a_j$ of $K$ is not isothetic to the vertex $a_j$ of $K'$. In this case \eqref{alt_uno} does not hold for all $u\in N(K,a_j)\cup N(K',a_j)$. If $a_j$ is a vertex of one polygon and it belongs to the relative interior of an edge of the other polygon, then the unit outer normal to this edge is necessarily  ${u^0_j}$. Since $\pa K$ and $\pa K'$ bifurcates at $a_j$, \eqref{alt_uno} does not hold when $u=u^0_j$.  

The identity \eqref{caratt_U} implies $u_1=u^0_1$ and $u_2=u^0_2$. The definition of $u_j^0$ clearly implies that $K_{u_j}$ and $K'_{u_j}$ are line segments, for $j=1,2$,  and that $K_{u_j}\cap K'_{u_j}$ is a line segment which, in counterclockwise order on $\pa K$, follows $a_j$ when $j=1$ and precedes $a_j$ when $j=2$. 
Assumption~\eqref{ass3}  of Lemma~\ref{archi_simmetrici} implies that the upper endpoints of $K_{u_1}$ and $K'_{u_1}$ coincide and that this point  belongs to $\relint\Si$. Similar arguments prove the analogous property for the lower endpoints of $K_{u_2}$ and $K'_{u_2}$. The properties proved up till now for $K$ and $K'$ can be proved, by similar arguments, also for $-L$ and $-L'$.

If, for some $j$, $\len(K_{u_j})\neq \len(K'_{u_j})$,  then \eqref{coni_giusti} is an immediate consequence of Alternative~\eqref{coincidono1} of Lemma~\ref{coincidono_intorno_lato}.
Assume $\len(K_{u_j})= \len(K'_{u_j})$. In this case the synisothesis of $(K,-L)$ and $(K',-L')$ implies $\len((-L)_{u_j})=\len((-L')_{u_j})$. 
Therefore, $a_j$ is a vertex of $K$ and of $K'$ and $b_j$ is a vertex of $-L$ and of $-L'$. Since $N(K,a_j)\cap S^1$, $N(K',a_j)\cap S^1$, $N(-L,b_j)\cap S^1$ and $N(-L',b_j)\cap S^1$ have the endpoint $u_j$ in common,  there exists a ``perturbation'' $\bar{u}$ of $u_j$ which belongs to the relative interior of $N(K,a_j)$, $N(K',a_j)$, $N(-L,b_j)$ and $N(-L',b_j)$.  The  vertex $a_j$ of $K$ is not isothetic to the vertex $a_j$ of $K'$, because $\pa K$ and  $\pa K'$ bifurcate at  $a_j$. 
Therefore,  the vertex $a_j$ of $K$ is isothetic to the vertex $b_j$ of $-L'$, and the vertex $a_j$ of $K'$ is isothetic to the vertex $b_j$ of $-L$. This is equivalent to  \eqref{coni_giusti}.  \end{proof}

Let $\mathcal{R}$ denote  clockwise rotation by $\pi/2$ and let
\begin{equation}\label{def_theta}
\te=\mathcal{R}\left( \frac{a_2-a_1}{\|a_2-a_1\|}\right)\quad\text{and}\quad
\te'=\mathcal{R}\left( \frac{b_2-b_1}{\|b_2-b_1\|}\right).
\end{equation}

Since the assumptions and the conclusion of Lemma~\ref{archi_simmetrici} are preserved by substituting  $K$, $-L$, $K'$ and $-L'$ with $-L$, $K$, $-L'$ and $K'$, respectively, we may assume $\te\leq\te'$ without loss of generality.
The next claim  states that, under suitable hypotheses, the arcs $\Si$ and $\Om+a_1-b_1$ can bifurcate only at  a point $c_1$ where every outer normal to these arcs  is $\te$, or it is larger than $\te'$.

\begin{claim}\label{claim2} Let $\te$ and $\te'$ be as in \eqref{def_theta}, with $\te\leq\te'$. For $i=1,2$, let  $a_i$, $a_i'$, $b_i$, $b_i'$ and $u_i$ be as in Claim~\ref{claim1} and let $x_i=b_i-a_i$. 
Assume that for some $j\in\{1,2\}$ the following two hypotheses hold:
if $\len(K_{u_j})\neq\len(K'_{u_j})$,  then Alternative~\eqref{coincidono1} of Lemma~\ref{coincidono_intorno_lato}, with $u$ replaced by $u_j$, holds;
we have $c_j\in\relint \Si\cap \relint (\Om+x_j)$, where $c_j$ is the endpoint different from $a_j$ of the arc that is contained in $\Si\cap (\Om+x_j)$ and contains $a_j$.
Then, when $j=1$ we have
\begin{align}
N(K,c_1)\cap N(-L+x_1,c_1)\cap[u_1,\te)_{S^1}&=\emptyset\quad\text{and}\label{intersezione}\\
N(K,c_1)\cap N(-L+x_1,c_1)\cap(\te,\te')_{S^1}&=\emptyset,\label{intersezione_due}
\end{align}
while when $j=2$ we have
\begin{align}
N(K,c_2)\cap N(-L+x_2,c_2)\cap(\te',u_2]_{S^1}&=\emptyset\quad\text{and}\label{intersezionep}\\
N(K,c_2)\cap N(-L+x_2,c_2)\cap(\te,\te')_{S^1}&=\emptyset.\label{intersezione_duep}
\end{align}
\end{claim}

\begin{proof}We prove the claim when $j=1$. To prove \eqref{intersezione} assume that there exists $w\in S^1$ in the intersection of the three sets. 
Let us first prove that $a_1$ is a vertex of $K$ and of $K'$ when $w=u_1$. Indeed, due to Claim~\ref{claim1}, $a_1$ is a vertex of both sets if and only if $\len(K_{u_1})=\len(K'_{u_1})$.   If $\len(K_{u_1})\neq\len(K'_{u_1})$,  then Alternative~\eqref{coincidono1} of Lemma~\ref{coincidono_intorno_lato} implies that $K$ and $-L'+x_1$ coincide in a neighbourhood of $K_{u_1}$. In particular, $c_1\notin[a_1,a_1']$. Thus $u_1\notin N(K,c_1)$, a contradiction to $w=u_1$. Similar arguments prove that  $b_1$ is a vertex of $-L$ and of $-L'$ when $w=u_1$.

Let $x_0=a_1+c_1-x_1$. We claim that
$g_{K,L}(x)\neq g_{K',L'}(x)$ for some $x$ close to $x_0$. First we prove \begin{equation}\label{identita_sets}
K\cap(L+x_0)=K'\cap(L'+x_0).
\end{equation}
The translation by $x_0$ maps the points $-c_1+x_1$ and $-b_1$ of $\pa L\cap \pa L'$ respectively to $a_1$ and $c_1$;  see Fig.~\ref{archi_claim2}. 
\begin{figure}
\begin{center}
\begin{picture}(0,0)%
\includegraphics{archi_claim2_2bis.pstex}%
\end{picture}%
\setlength{\unitlength}{3522sp}%
\begingroup\makeatletter\ifx\SetFigFontNFSS\undefined%
\gdef\SetFigFontNFSS#1#2#3#4#5{%
  \reset@font\fontsize{#1}{#2pt}%
  \fontfamily{#3}\fontseries{#4}\fontshape{#5}%
  \selectfont}%
\fi\endgroup%
\begin{picture}(4227,1858)(-2443,766)
\put(-1942,1175){\makebox(0,0)[lb]{\smash{{\SetFigFontNFSS{10}{12.0}{\familydefault}{\mddefault}{\updefault}$w$}}}}
\put(-2428,1636){\makebox(0,0)[lb]{\smash{{\SetFigFontNFSS{10}{12.0}{\familydefault}{\mddefault}{\updefault}$u_1$}}}}
\put(-1651,1225){\makebox(0,0)[lb]{\smash{{\SetFigFontNFSS{10}{12.0}{\familydefault}{\mddefault}{\updefault}$\theta$}}}}
\put(954,2012){\makebox(0,0)[lb]{\smash{{\SetFigFontNFSS{10}{12.0}{\familydefault}{\mddefault}{\updefault}$a_2$}}}}
\put(-257,2429){\makebox(0,0)[lb]{\smash{{\SetFigFontNFSS{12}{14.4}{\familydefault}{\mddefault}{\updefault}$\pa K$}}}}
\put(1455,2009){\makebox(0,0)[lb]{\smash{{\SetFigFontNFSS{12}{14.4}{\familydefault}{\mddefault}{\updefault}$\pi$}}}}
\put(1009,1312){\makebox(0,0)[lb]{\smash{{\SetFigFontNFSS{12}{14.4}{\familydefault}{\mddefault}{\updefault}$-\pi+x_1+x_0$}}}}
\put(-56,848){\makebox(0,0)[lb]{\smash{{\SetFigFontNFSS{12}{14.4}{\familydefault}{\mddefault}{\updefault}$\pa L'+x_0$}}}}
\put(-701,2429){\makebox(0,0)[lb]{\smash{{\SetFigFontNFSS{12}{14.4}{\familydefault}{\mddefault}{\updefault}$\pa K'$}}}}
\put(720,1029){\makebox(0,0)[lb]{\smash{{\SetFigFontNFSS{12}{14.4}{\familydefault}{\mddefault}{\updefault}$\pa L+x_0$}}}}
\put(-518,1986){\makebox(0,0)[lb]{\smash{{\SetFigFontNFSS{10}{12.0}{\familydefault}{\mddefault}{\updefault}$a_1$}}}}
\put(347,1322){\makebox(0,0)[lb]{\smash{{\SetFigFontNFSS{10}{12.0}{\familydefault}{\mddefault}{\updefault}$c_1$}}}}
\end{picture}%
\end{center}
\caption{$K\cap(L+x_0)$ (dotted lines) and $K'\cap (L'+x_0)$ (continuous lines).}
\label{archi_claim2}
\end{figure}
Let $\pi=\set{y\in \Real^2 :  (y-a_1)\cdot w \geq 0 }$. 
The sets $K$, $K'$, $-L+x_1$ and $-L'+x_1$ are contained in $-\pi+x_1+x_0=\{y\in\Real^2 : (y-c_1)\cdot w\leq 0\}$ (because $w$ is an outer normal to all these sets at $c_1$). 
The inclusions  $-L+x_1$, $-L'+x_1\subset-\pi+x_1+x_0$ are equivalent to $L+x_0$, $L'+x_0\subset\pi$. Therefore
\begin{equation}\label{incluso_strip}
 K\cap(L+x_0), K'\cap(L'+x_0)\subset\pi\cap(-\pi+x_1+x_0).
\end{equation}
We claim that
\begin{equation}\label{diff_simm1}
K\cap\pi=K'\cap\pi\quad\text{and }\quad\cl(K\De K')\cap\pi=\{a_1\}.
\end{equation}
When $w=u_1$ the first identity is true because $K\cap\pi=[a_1,a_1']=K'\cap\pi$, by Claim~\ref{claim1} and because  $a_1$ is a vertex of $K$ and of $K'$.  When $w\in(u_1,\te)_{S^1}$,  we have $\pa K\cap\pi=\pa K'\cap\pi$, because $\pa K\cap\pi$ and $\pa K'\cap\pi$ are contained in $\Si$ (observe that $a_2\notin\pi$ because $w<\te$), which is contained in $\pa K\cap \pa K'$. Since $K$ and $K'$ are convex, the identity $\pa K\cap\pi=\pa K'\cap\pi$ implies $K\cap\pi=K'\cap\pi$. To prove the second identity observe that $K\De K'$ is contained in one of the halfplanes bounded by the line through $a_1$ and $a_2$, and observe that this halfplane intersects $\pi\cap K$ only in $a_1$.
Similar arguments prove $(-L+x_1)\cap\pi=(-L'+x_1)\cap\pi$ and $\cl((-L+x_1)\De (-L'+x_1))\cap\pi=\{a_1\}$. These identities are equivalent to 
\begin{gather}
(L+x_0)\cap(-\pi+x_1+x_0)=(L'+x_0)\cap(-\pi+x_1+x_0)\quad\text{and}\label{diff_simm2}\\
\cl((L+x_0)\De (L'+x_0))\cap(-\pi+x_1+x_0)=\{c_1\}.\label{diff_simm3}
\end{gather}
Formulas  \eqref{incluso_strip}, \eqref{diff_simm1} and \eqref{diff_simm2} imply \eqref{identita_sets}.

Formulas \eqref{diff_simm1} and \eqref{diff_simm3} imply that $\cl( (K\Delta K')\cup ((L+x_0)\Delta (L'+x_0)))$ intersects the strip $\pi\cap(-\pi+x_1+x_0)$ only in $a_1$ and $c_1$. Therefore, when $x$ is close to $x_0$, the sets $K\cap(L+x)$ and $K'\cap(L'+x)$ may differ only in a neighbourhood of $a_1$ and of $c_1$ (see Fig.~\ref{archi_claim2_3}). 
\begin{figure}
\begin{center}
\begin{picture}(0,0)%
\includegraphics{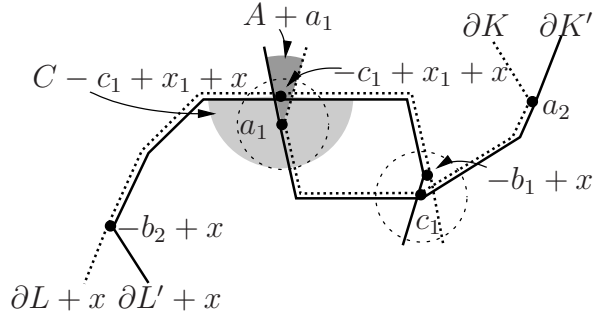}%
\end{picture}%
\setlength{\unitlength}{4144sp}%
\begingroup\makeatletter\ifx\SetFigFontNFSS\undefined%
\gdef\SetFigFontNFSS#1#2#3#4#5{%
  \reset@font\fontsize{#1}{#2pt}%
  \fontfamily{#3}\fontseries{#4}\fontshape{#5}%
  \selectfont}%
\fi\endgroup%
\begin{picture}(3298,1921)(-1721,-2129)
\put(-348,-1014){\makebox(0,0)[lb]{\smash{{\SetFigFontNFSS{12}{14.4}{\familydefault}{\mddefault}{\updefault}$a_1$}}}}
\put(693,-1582){\makebox(0,0)[lb]{\smash{{\SetFigFontNFSS{12}{14.4}{\familydefault}{\mddefault}{\updefault}$c_1$}}}}
\put(1442,-898){\makebox(0,0)[lb]{\smash{{\SetFigFontNFSS{12}{14.4}{\familydefault}{\mddefault}{\updefault}$a_2$}}}}
\put(1425,-429){\makebox(0,0)[lb]{\smash{{\SetFigFontNFSS{12}{14.4}{\familydefault}{\mddefault}{\updefault}$\pa K'$}}}}
\put(-329,-367){\makebox(0,0)[lb]{\smash{{\SetFigFontNFSS{12}{14.4}{\familydefault}{\mddefault}{\updefault}$A+a_1$}}}}
\put(970,-429){\makebox(0,0)[lb]{\smash{{\SetFigFontNFSS{12}{14.4}{\familydefault}{\mddefault}{\updefault}$\pa K$}}}}
\put(199,-699){\makebox(0,0)[lb]{\smash{{\SetFigFontNFSS{12}{14.4}{\familydefault}{\mddefault}{\updefault}$-c_1+x_1+x$}}}}
\put(-1078,-1632){\makebox(0,0)[lb]{\smash{{\SetFigFontNFSS{12}{14.4}{\familydefault}{\mddefault}{\updefault}$-b_2+x$}}}}
\put(1110,-1344){\makebox(0,0)[lb]{\smash{{\SetFigFontNFSS{12}{14.4}{\familydefault}{\mddefault}{\updefault}$-b_1+x$}}}}
\put(-1574,-736){\makebox(0,0)[lb]{\smash{{\SetFigFontNFSS{12}{14.4}{\familydefault}{\mddefault}{\updefault}$C-c_1+x_1+x$}}}}
\put(-1706,-2051){\makebox(0,0)[lb]{\smash{{\SetFigFontNFSS{12}{14.4}{\familydefault}{\mddefault}{\updefault}$\pa L+x$}}}}
\put(-1069,-2056){\makebox(0,0)[lb]{\smash{{\SetFigFontNFSS{12}{14.4}{\familydefault}{\mddefault}{\updefault}$\pa L'+x$}}}}
\end{picture}%
\end{center}
\caption{$K\cap(L+x)$ (dotted lines) and $K'\cap (L'+x)$ (continuous lines) for a suitable $x$ close to $x_0$.}
\label{archi_claim2_3}
\end{figure}
Assume $K\subset K'$ in a neighbourhood of $a_1$ (note that either we have $K\subset K'$ or we have $K'\subset K$, because $\pa K$ and $\pa K'$ coincide on one side of $a_1$). Define $A=\cone(K',a_1)\setminus\cone(K,a_1)$. By \eqref{coni_giusti} we have $A=\cone(-L,b_1)\setminus \cone(-L',b_1)$.
Let $C=\cone(L,x_1-c_1)=\cone(L',x_1-c_1)$ and  $-D=\cone(K,c_1)=\cone(K',c_1)$ (recall that $K=K'$ near $c_1$ and   $L=L'$ near $x_1-c_1$).
If $\de>0$ is chosen sufficiently small, outside $B(a_1,\de)\cup B(c_1,\de)$ we have $K\cap(L+x)=K'\cap(L'+x)$. 
On the other hand
\begin{equation*}
\left((K'\cap(L'+x))\setminus(K\cap(L+x))\right)\cap B(a_1,\de)=
(A+a_1)\cap (C+x_1-c_1+x);
\end{equation*}
see Fig.~\ref{archi_claim2_3}. Therefore 
\begin{align*}
\area(K\cap(L+x)\cap B(a_1,\de))-\area&(K'\cap(L'+x)\cap B(a_1,\de))\\
&=-\area((A+a_1)\cap(C+x_1-c_1+x))\\
&=-\area(A\cap(C+x-x_0))\\
&=-g_{A,C}(x-x_0).
\end{align*}
Similar arguments prove
\begin{equation*}
\area(K\cap(L+x)\cap B(c_1,\de))-\area(K'\cap(L'+x)\cap B(c_1,\de))=g_{A,D}(x-x_0).
\end{equation*}
Therefore, for each $x$ in a neighbourhood of $x_0$, we have
\begin{equation}\label{conn_cov_coni}
g_{K,L}(x)-g_{K',L'}(x)=g_{A,D}(x-x_0)-g_{A,C}(x-x_0).
\end{equation}
We apply Lemma~\ref{coni_finale}, with $B=\{O\}$, $n=2$ and the Cartesian coordinates chosen so that $w=(0,-1)$. The assumptions of this lemma are satisfied. Indeed, $C,D\subset \{x_2\leq0\}$, and either $C\subset D$ or $D\subset C$, because the lower endpoints of $C\cap S^1$ and of $D\cap S^1$ coincide, by $[a_1,c_1]_{\pa K}\subset\pa K\cap\pa (-L+x_1)$.  Moreover, $\cl(K\Delta K')\cap\pi=\{a_1\}$ implies $A\cap \{x_2=0\}=\{O\}$.
Observe that we have $C\neq D$, because $\pa K$ and $\pa(-L)+x_1$ bifurcate at $c_1$. Thus, this lemma implies $g_{A,D}\not\equiv g_{A,C}$. Since $g_{A,D}$ and $g_{A,C}$ are homogeneous functions of degree $2$, they do not coincide in any neighbourhood of $O$.  Thus \eqref{conn_cov_coni} implies $g_{K,L}\neq g_{K',L'}$. This contradiction proves \eqref{intersezione}.

The proof of \eqref{intersezione_due} is similar although simpler.  Assume that $w\in S^1$ belongs to the three sets in \eqref{intersezione_due} and define $x_0=a_2+c_1-x_1$. We prove again that $g_{K,L}$ does not coincide with $g_{K',L'}$ in a neighbourhood of $x_0$. The translation by $x_0$ maps the point $-c_1+x_1$ of $\pa L\cap\pa L'$  to the point  $a_2$. The identity $K\cap(L+x_0)=K'\cap(L'+x_0)$ is proved as before.  Defining $\pi=\set{y : (y-a_2)\cdot w \geq 0}$, the inclusion \eqref{incluso_strip} holds also in this case.
Moreover, when $x$ is close to $x_0$ the sets $K\cap(L+x)$ and $K'\cap(L'+x)$ may differ only in a neighbourhood of $a_2$, because $\cl(K\Delta K')$ intersects the strip  $\pi\cap (-\pi+x_1+x_0)$ only in $a_2$ (recall that $a_1\notin\pi$ because $w>\te$) and $\cl ((L+x_0)\Delta (L'+x_0))$ does not intersect this strip (recall that $-b_1+x_0\notin -\pi+x_1+x_0$ because this is equivalent to $a_1\notin\pi$ which holds true, and $-b_2+x_0\notin -\pi+x_1+x_0$ because $w<\te'$).
Assume $K\subset K'$ in a neighbourhood of $a_2$, let $A=\cone(K',a_2)\setminus\cone(K,a_2)$ and
let $C=\cone(L,x_1-c_1)=\cone(L',x_1-c_1)$.  
For each $x$ close to $x_0$ we have
\[
g_{K,L}(x)-g_{K',L'}(x)=-g_{A,C}(x-x_0).
\]
Lemma~\ref{coni_finale}, with $D=B=\set{O}$, implies that the previous formula contradicts $g_{K,L}=g_{K',L'}$. This contradiction proves \eqref{intersezione_due}.
\end{proof}

\begin{claim}\label{claim2a} Let $\cl U=[u_1,u_2]$.
If  $\len(K_{u_j})\neq\len(K'_{u_j})$, for some $j\in\{1,2\}$, then Alternative~\eqref{coincidono1} of Lemma~\ref{coincidono_intorno_lato}, with $u$ substituted by $u_j$, holds.
\end{claim}

\begin{proof}Let $\te$ and $\te'$ be as in \eqref{def_theta}. We may assume $\te\leq\te'$. For $i=1,2$, let  $a_i$, $a_i'$, $b_i$ and $b_i'$ be as in Claim~\ref{claim1} and let $x_i=b_i-a_i$.  Assume the claim false for some $j$. 
By \eqref{alt_due}, with $u$ substituted by $u_j$, we have 
\begin{equation}\label{lunghezze}
\len(K_{u_j})=\len((-L')_{u_j})\quad\text{and}\quad\len(K'_{u_j})=\len((-L)_{u_j}).
\end{equation}
Let $a_j''$ be the endpoint of $K_{u_j}$ which does not belong to $K'_{u_j}$ or the endpoint of $K'_{u_j}$ which does not belong to $K_{u_j}$, according to whether $\len(K_{u_j})>\len(K'_{u_j})$ or $\len(K_{u_j})<\len(K'_{u_j})$. Define  $b_j''$ similarly, by substituting $K$ with $-L$ and $K'$ with $-L'$. 

First we prove that if Claim~\ref{claim3} is false both when $j=1$ and when $j=2$, then  $(K,-L)$ and $(K',-L')$ are not synisothetic, contrary to what assumed in Lemma~\ref{archi_simmetrici}. 
Let $w\in S^1$ be orthogonal to the edges of $K$ adjacent to $K_{u_2}$ and to those of $K'$ adjacent  to $K'_{u_2}$, oriented in such a way that $w\in N(K,a'_2)$. We recall that Lemma~\ref{coincidono_intorno_lato} proves that those edges are parallel.
Define $v$ as $w$, with  $K$, $K'$ and $a'_2$ replaced by  $-L$, $-L'$ and $b'_2$, respectively. We stress that $w\neq v$, because otherwise both \eqref{coincidono1} and \eqref{coincidono2} in  Lemma~\ref{coincidono_intorno_lato} hold true, and this implies that Claim~\ref{claim3} is true when $j=2$. We distinguish two possible cases.

\textit{Case $a_1'=a_2'$.} We have $\Si=[a_1,a_1']\cup[a_2',a_2]$, $w=u_1$ and $u_2\in(u_1,-u_1)_{S^1}$; see Fig.~\ref{fig_claim2a}. 

\textit{Case $a_1'\neq a_2'$.} We have $\Si=[a_1,a_1']\cup[a_1',a_2']\cup[a_2',a_2]$, $w\neq u_1$ and $w\neq u_2$. Moreover, the parallelism of the edges of $K$ adjacent to $K_{u_i}$ and the parallelism of the edges of $K'$ adjacent to $K'_{u_i}$, for each $i\in\{1,2\}$, together with the convexity of $K$ and $K'$, imply that $K=\conv(a_1,a_1',a_2',a_2)$ and $K'=\conv(a_1'',a_1',a_2',a_2'')$, or vice versa.

Similar descriptions hold also for $-L$ and $-L'$, with  $v$, $b_i$, $b_i'$ and $b_i''$ replacing, respectively, $w$, $a_i$, $a_i'$ and $a_i''$.

\begin{figure}
\begin{center}
\begin{picture}(0,0)%
\includegraphics{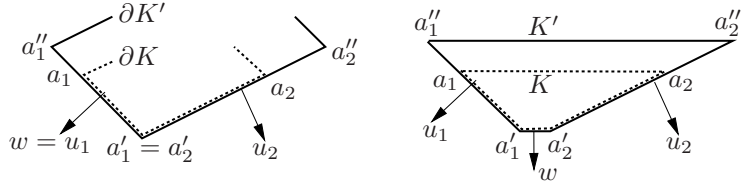}%
\end{picture}%
\setlength{\unitlength}{3315sp}%
\begingroup\makeatletter\ifx\SetFigFontNFSS\undefined%
\gdef\SetFigFontNFSS#1#2#3#4#5{%
  \reset@font\fontsize{#1}{#2pt}%
  \fontfamily{#3}\fontseries{#4}\fontshape{#5}%
  \selectfont}%
\fi\endgroup%
\begin{picture}(5392,1405)(1246,-665)
\put(4857,-402){\makebox(0,0)[lb]{\smash{{\SetFigFontNFSS{10}{12.0}{\familydefault}{\mddefault}{\updefault}$a_1'$}}}}
\put(5246,-397){\makebox(0,0)[lb]{\smash{{\SetFigFontNFSS{10}{12.0}{\familydefault}{\mddefault}{\updefault}$a_2'$}}}}
\put(4397, 92){\makebox(0,0)[lb]{\smash{{\SetFigFontNFSS{10}{12.0}{\familydefault}{\mddefault}{\updefault}$a_1$}}}}
\put(6139, 86){\makebox(0,0)[lb]{\smash{{\SetFigFontNFSS{10}{12.0}{\familydefault}{\mddefault}{\updefault}$a_2$}}}}
\put(6479,518){\makebox(0,0)[lb]{\smash{{\SetFigFontNFSS{10}{12.0}{\familydefault}{\mddefault}{\updefault}$a_2''$}}}}
\put(6121,-376){\makebox(0,0)[lb]{\smash{{\SetFigFontNFSS{10}{12.0}{\familydefault}{\mddefault}{\updefault}$u_2$}}}}
\put(4321,-286){\makebox(0,0)[lb]{\smash{{\SetFigFontNFSS{10}{12.0}{\familydefault}{\mddefault}{\updefault}$u_1$}}}}
\put(5176,-601){\makebox(0,0)[lb]{\smash{{\SetFigFontNFSS{10}{12.0}{\familydefault}{\mddefault}{\updefault}$w$}}}}
\put(4258,516){\makebox(0,0)[lb]{\smash{{\SetFigFontNFSS{10}{12.0}{\familydefault}{\mddefault}{\updefault}$a_1''$}}}}
\put(1353,337){\makebox(0,0)[lb]{\smash{{\SetFigFontNFSS{10}{12.0}{\familydefault}{\mddefault}{\updefault}$a_1''$}}}}
\put(1541,121){\makebox(0,0)[lb]{\smash{{\SetFigFontNFSS{10}{12.0}{\familydefault}{\mddefault}{\updefault}$a_1$}}}}
\put(3061,-421){\makebox(0,0)[lb]{\smash{{\SetFigFontNFSS{10}{12.0}{\familydefault}{\mddefault}{\updefault}$u_2$}}}}
\put(3196, 29){\makebox(0,0)[lb]{\smash{{\SetFigFontNFSS{10}{12.0}{\familydefault}{\mddefault}{\updefault}$a_2$}}}}
\put(3646,299){\makebox(0,0)[lb]{\smash{{\SetFigFontNFSS{10}{12.0}{\familydefault}{\mddefault}{\updefault}$a_2''$}}}}
\put(1261,-376){\makebox(0,0)[lb]{\smash{{\SetFigFontNFSS{10}{12.0}{\familydefault}{\mddefault}{\updefault}$w=u_1$}}}}
\put(1981,-421){\makebox(0,0)[lb]{\smash{{\SetFigFontNFSS{10}{12.0}{\familydefault}{\mddefault}{\updefault}$a_1'=a_2'$}}}}
\put(2071,569){\makebox(0,0)[lb]{\smash{{\SetFigFontNFSS{10}{12.0}{\familydefault}{\mddefault}{\updefault}$\pa K'$}}}}
\put(2071,254){\makebox(0,0)[lb]{\smash{{\SetFigFontNFSS{10}{12.0}{\familydefault}{\mddefault}{\updefault}$\pa K$}}}}
\put(5096,474){\makebox(0,0)[lb]{\smash{{\SetFigFontNFSS{10}{12.0}{\familydefault}{\mddefault}{\updefault}$K'$}}}}
\put(5101, 44){\makebox(0,0)[lb]{\smash{{\SetFigFontNFSS{10}{12.0}{\familydefault}{\mddefault}{\updefault}$K$}}}}
\end{picture}%
\end{center}
\caption{The two possible configurations for $K$ and $K'$.}
\label{fig_claim2a}
\end{figure}

It cannot be $a_1'=a_2'$ and $b_1'=b_2'$, because the descriptions above imply $w=v$, which is false.
Assume $a_1'\neq a_2'$ and $b_1'\neq b_2'$. If, say, $K'=\conv(a_1'',a_1',a_2',a_2'')$, then no edge of  $-L$  is a translate of the edge $[a_1'',a_2'']$ of $K'$, because no edge of $-L$ is orthogonal to $w$. Thus, the synisothesis of $(K,-L)$ and $(K',-L')$ implies that an edge of $K$ is a translate of $[a_1'',a_2'']$. This implies that  $u_2=-u_1$. Both $K$, $K'$, $L$ and  $L'$ are parallelograms with two edges orthogonal to $u_1$. This property and  \eqref{lunghezze} imply that, up to an affine transformation,  $(K,L)$ and $(K',L')$ are trivial associates of $(\cK_3,\cL_3)$ and $(\cK_4,\cL_4)$, respectively.
This contradicts the assumptions of Lemma~\ref{archi_simmetrici}.

Assume $a_1'\neq a_2'$ and $b_1'=b_2'$. We have $\len([a_1'',a_2''])>\len([a_1,a_2])>\len([a_1',a_2'])$, because $u_2\in(u_1,-u_1)_{S^1}$.  We also have $\len([a_1'',a_2''])=\len([b_1'',b_2''])+\len([a_1',a_2'])>\len([b_1'',b_2''])>\len([b_1',b_2'])$, by the previous descriptions, \eqref{lunghezze} with $j=1$ and \eqref{lunghezze} with $j=2$. Therefore, if, say,  $K'=\conv(a_1'',a_1',a_2',a_2'')$, then neither $K$ nor $-L$ have an edge which is a translate of the edge $[a_1'',a_2'']$ of $K'$. This contradicts the synisothesis of $(K,-L)$ and $(K',-L')$.

It remains to prove that it is not possible that Claim~\ref{claim3} holds for one index, say $j=1$, and it does not hold for the other one, say $j=2$. Assume this false and let $w$ and $v$ be defined as above. Up to exchanging the roles of $K$ and $K'$ we may assume $\len(K_{u_2})<\len(K'_{u_2})$, so that $a_2$ is a vertex of $K$ and $N(K,a_2)\cap S^1=[u_2,-w]_{S^1}$. Therefore, since $[a_1,a_2]$ is a chord of $K$, we have 
\begin{equation}\label{wminte}
\text{$w\leq\te$ and $w=\te$ if and only if $[a_1,a_2]$ is an edge of $K$.}
\end{equation}
Similar arguments prove that 
\begin{equation}\label{vminte}
\text{$v\leq\te'$ and $v=\te'$ if and only if $[b_1,b_2]$ is an edge of $-L'$.}
\end{equation}
We claim that $a_2'\notin [a_1,c_1)_{\pa K}$. If $a_2'\in [a_1,c_1)_{\pa K}$, then $\pa K$ and $\pa(-L')+x_1$ coincide in a neighbourhood of $a_2'$. In particular $-L'$ has a vertex $p$ with $N(-L',p)\cap S^1=[w,u_2]_{S^1}$. This is false, because $b_2'$ is the only vertex of $-L'$ with the property that $u_2$ is the upper endpoint of the intersection of $S^1$ with the normal cone at that vertex, however the lower endpoint of $N(-L,b_2')$ is $v$ and  $w\neq v$. This contradiction proves $a_2'\notin [a_1,c_1)_{\pa K}$.

Let $z$ be the lower endpoint of $N(K,c_1)\cap S^1$. The vector $z$ is also the lower endpoint of $N(-L'+x_1,c_1)\cap S^1$, since $[a_1,c_1]_{\pa K}\subset\pa K\cap (\pa(-L')+x_1)$. Since $a_2'\notin [a_1,c_1)_{\pa K}$ we have $z\leq w$, by convexity.  Thus, \eqref{intersezione} and \eqref{wminte} imply $z=w=\te$.
The line segment $[a_1,a_2]$ is an edge of $K$, by \eqref{wminte}. Therefore, by \eqref{coni_giusti} with $j=1$,  the line $l$ through $b_1$ and orthogonal to $\te$ supports $-L'$. This implies $\te\geq\te'$. Since we assumed $\te\leq\te'$, we have $\te=\te'$. Thus $b_2\in l$ and  $[b_1,b_2]$ is an edge of $-L'$. The equalities $\te'=v$ (a consequence of \eqref{vminte}), $\te=w$ and $\te=\te'$ contradict $v\neq w$. \end{proof}

Claims~\ref{claim1} and~\ref{claim2a} imply that \eqref{coni_giusti} holds both when $j=1$ and when $j=2$. 

\begin{claim}\label{claim3}Let $\cl U=[u_1,u_2]$ and, for $i=1,2$, let $x_i=a_i-b_i$, where $a_i$ and $b_i$ are as in Claim~\ref{claim1}. 
Assume that $\Si$ is not a translate of $\Om$. Then $U$ contains an half-circle and $x_1\neq x_2$. 
Moreover, if $U$ is an half-circle, then $\Si$ is the union of three consecutive edges $L_1$, $M$ and $L_2$ of a parallelogram, $L_1$ and $L_2$ are orthogonal to $u_1$ and $\Om$ is the union of $L_1-x_1$, of a line segment $M'$ parallel to $M$ and of $L_2-x_2$.\end{claim}

\begin{proof} Let $\te$ and $\te'$ be as in \eqref{def_theta} and, for $i=1,2$, let  $a_i'$ and $b_i'$ be as in Claim~\ref{claim1} and let $c_i$ be as in Claim~\ref{claim2}. We may assume $\te\leq\te'$. First we prove 
\begin{equation}\label{ai_in_relint}
c_i\in\relint\Si\cap(\relint\Om+x_i),
\end{equation}
for $i=1,2$. Assume \eqref{ai_in_relint} false when $i=2$. Since $\Si$ and $\Om+x_2$ coincide in a neighbourhood of their common upper endpoint $a_2$, by Claim~\ref{claim1},  $c_2$ coincides with the lower endpoint of one arc, that is,  with $a_1$ or with $b_1+x_2$. 
Assume $c_2=a_1$. In this case $\Si\subset\Om+x_2$ and, since $\Si\neq\Om+x_2$ by assumption, $a_1\neq b_1+x_2$. 
The portion of $\Om$ with  outer normal $u_1$, that is $[a_1,a_1']$, should be contained in the portion of $\Om+x_2$ with outer normal $u_1$, that is in $[b_1,b_1']+x_2$. Moreover, the inequality $a_1\neq b_1+x_2$ implies $\len([a_1,a_1'])<\len([b_1,b_1'])$.
The previous description implies that $\Om+x_1$  bifurcates from $\Si$ at $a_1'$, that is, it implies $c_1=a_1'$. Claim~\ref{claim2a} and the observation that $c_1=a_1'\in\relint\Si\cap(\relint\Om+x_1)$ imply that \eqref{intersezione} holds. On the other hand, $c_1=a_1'$  implies  $u_1\in N(K,c_1)\cap N(-L+x_1,c_1)$, which contradicts \eqref{intersezione}. Similar arguments prove $c_2\neq b_1+x_2$, and prove \eqref{ai_in_relint} when $i=1$.

Let $w_1$ be the lower endpoint of $N(K,c_1)\cap S^1$ and of $N(-L'+x_1,c_1)\cap S^1$ (these endpoints coincide because $[a_1,c_1]_{\pa K}\subset \Om\cap(\pa(-L)+x_1)$). Let $w_2$ be the upper endpoint of $N(K,c_2)\cap S^1$ and $N(-L'+x_2,c_2)\cap S^1$. We prove
\begin{equation}\label{w1minw2}
 w_1\leq w_2.
\end{equation}
Assume  $w_1>w_2$. In this case the two sub-arcs $[a_1,c_1]_{\pa K}$ and $[c_2,a_2]_{\pa K}$ of $\pa K$  overlap and contain  the arc $[c_2,c_1]_{\pa K}$. Moreover $[c_2,c_1]_{\pa K}$ is neither a point nor a line segment, because it contains a segment orthogonal to $w_i$, for $i=1,2$, by definition of $w_i$. Therefore the inclusions
\begin{equation*}
[c_2,c_1]_{\pa K}\subset[a_1,c_1]_{\pa K}\subset\Om+x_1\quad \text{and}\quad
[c_2,c_1]_{\pa K}\subset[c_2,a_2]_{\pa K}\subset\Om+x_2,
\end{equation*}
imply $x_1=x_2$. Since $\Om+x_1$ is contained in  the union of two overlapping sub-arcs of $\Si$, and its endpoints coincide with those of $\Sigma$, we have $\Om+x_1=\Si$. This contradicts the assumptions of Claim~\ref{claim3} and proves \eqref{w1minw2}.

Claim~\ref{claim2a} and~\eqref{ai_in_relint} imply that the assumptions of Claim~\ref{claim2} are satisfied for each $j=1,2$. Thus \eqref{intersezione}, \eqref{intersezione_due}, \eqref{intersezionep} and \eqref{intersezione_duep} hold. 

Assume  $\te=\te'$. Formulas~\eqref{intersezione} and~\eqref{intersezionep}  imply respectively  $w_1\geq \te$ and $w_2\leq\te$. Thus \eqref{w1minw2} implies  $w_1=\te=w_2$. 
Let  $K_\te=[\dot{a_1},\dot{a_2}]$ and $(-L)_\te=[\dot{b_1},\dot{b_2}]$, where $\dot{a_2}$ follows $\dot{a_1}$ and $\dot{b_2}$ follows $\dot{b_1}$, in counterclockwise order on the respective boundaries. We may clearly write
\begin{equation}\label{rappr_new}\begin{aligned}
\Si&=[a_1,\dot{a_1}]_{\pa K}\cup[\dot{a_1},\dot{a_2}]\cup[\dot{a_2},a_2]_{\pa K}\quad\text{and}\\
\Om&=[b_1,\dot{b_1}]_{\pa (-L)}\cup[\dot{b_1},\dot{b_2}]\cup[\dot{b_2},b_2]_{\pa (-L)},
\end{aligned}\end{equation}
where $[\dot{a_1},\dot{a_2}]$ and $[\dot{b_1},\dot{b_2}]$ are parallel.
The definitions of  $w_1$ and $w_2$ imply $c_i\in [\dot{a_1},\dot{a_2}]$ and $c_i\in[\dot{b_1},\dot{b_2}]+x_i$, for $i=1,2$, $c_1\neq \dot{a_1}$ and $c_2\neq \dot{a_2}$. Therefore, by definition of $c_1$ and of $c_2$ we have
\begin{equation}\label{rappr_new2}
[b_1,\dot{b_1}]_{\pa(-L)}=[a_1,\dot{a_1}]_{\pa K}-x_1\quad\text{and}\quad
[\dot{b_2},b_2]_{\pa (-L)}=[\dot{a_2},a_2]_{\pa K}-x_2.
\end{equation}
We have $\len([\dot{a_1},\dot{a_2}])\neq\len([\dot{b_1},\dot{b_2}])$, because otherwise $\Si=\Om+x_1$, by \eqref{rappr_new} and \eqref{rappr_new2}.

In order to prove that $U$ contains an half-circle it suffices to prove 
\begin{equation}\label{base_piu_lunga}
\|a_2-a_1\|\leq\|\dot{a_2}-\dot{a_1}\|.
\end{equation}
Indeed, if $U$ does not contain an half-circle, then the line $l$ through $a_1$ orthogonal to $u_1$ and the line $r$ through $a_2$ orthogonal to $u_2$ bound a cone with apex contained in the halfplane bounded by the line through $a_1$ and $a_2$ and containing $[\dot{a_1},\dot{a_2}]$. Since $l$ and $r$ support $K$, the cone contains  $[\dot{a_1},\dot{a_2}]$ and \eqref{base_piu_lunga} is false.

We assume $\|a_2-a_1\|>\|\dot{a_2}-\dot{a_1}\|$ and obtain a contradiction by proving that $g_{K,L}(x)\neq g_{K',L'}(x)$, for some $x$ close to $x_0$, where $x_0$ will be defined later. We may write
\begin{equation*}
a_2-a_1=\al(\dot{a_2}-\dot{a_1})\quad\text{ and }\quad \dot{b_2}-\dot{b_1}=\be(\dot{a_2}-\dot{a_1}),
\end{equation*}
with $\al>1$, $\be>0$. The inequality $\len([\dot{a_1},\dot{a_2}])\neq\len([\dot{b_1},\dot{b_2}])$ implies $\be\neq1$. Up to exchanging $K$, $K'$, $a_i$ and $\dot{a_i}$ with $-L$, $-L'$, $b_i$ and $\dot{b_i}$, respectively, we may assume $\be<1$.
Let $x_0=\dot{a_1}+b_1+\ga(\dot{a_2}-\dot{a_1})$, where $\ga$ is a fixed number in $(\be,\min(1,\be+\al-1))$. Easy computations  and the equalities $a_2-\dot{a_2}=b_2-\dot{b_2}$ and $a_1-\dot{a_1}=b_1-\dot{b_1}$ (that are consequences of \eqref{rappr_new2}) give the following expressions:
$-b_1+x_0=\dot{a_1}+\ga(\dot{a_2}-\dot{a_1})$; $-\dot{b_1}+x_0=a_1+(\ga/\al)(a_2-a_1)$; $-\dot{b_2}+x_0=a_1+((\ga-\be)/\al)({a_2}-{a_1})$ and
$-b_2+x_0=\dot{a_1}+(1-\al+\ga-\be)(\dot{a_2}-\dot{a_1})$.
Since $0<(\ga-\be)/\al<\ga/\al<1$ and $1-\al+\ga-\be<0$, the previous expressions imply the following formulas:
\begin{gather}
-\dot{b_1}+x_0,-\dot{b_2}+x_0\in\relint [a_1,a_2];\label{sit_x0_1} \\
-b_2+x_0\notin K\cup K';\label{sit_x0_2}\\
-b_1+x_0\in \relint [\dot{a_1},\dot{a_2}];\label{sit_x0_3}\\ \|-b_1+x_0-\dot{a_1}\|>\|\dot{b_2}-\dot{b_1}\|.\label{sit_x0_4}
\end{gather}
\begin{figure}
\begin{center}
\begin{picture}(0,0)%
\includegraphics{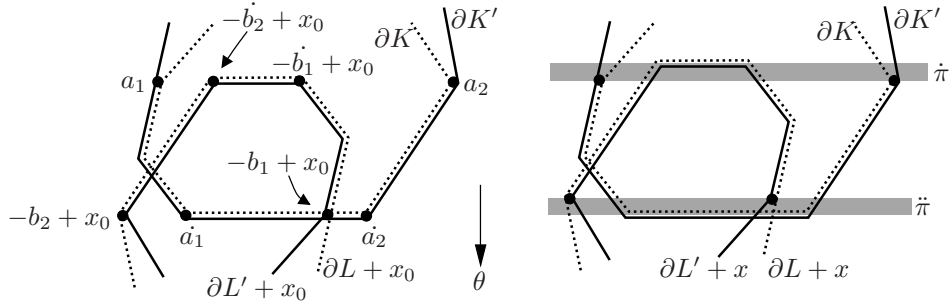}%
\end{picture}%
\setlength{\unitlength}{4144sp}%
\begingroup\makeatletter\ifx\SetFigFontNFSS\undefined%
\gdef\SetFigFontNFSS#1#2#3#4#5{%
  \reset@font\fontsize{#1}{#2pt}%
  \fontfamily{#3}\fontseries{#4}\fontshape{#5}%
  \selectfont}%
\fi\endgroup%
\begin{picture}(5510,1816)(-3309,-2050)
\put(1951,-386){\makebox(0,0)[lb]{\smash{{\SetFigFontNFSS{10}{12.0}{\familydefault}{\mddefault}{\updefault}$\pa K'$}}}}
\put(1496,-456){\makebox(0,0)[lb]{\smash{{\SetFigFontNFSS{10}{12.0}{\familydefault}{\mddefault}{\updefault}$\pa K$}}}}
\put(566,-1871){\makebox(0,0)[lb]{\smash{{\SetFigFontNFSS{10}{12.0}{\familydefault}{\mddefault}{\updefault}$\pa L'+x$}}}}
\put(2076,-1501){\makebox(0,0)[lb]{\smash{{\SetFigFontNFSS{10}{12.0}{\familydefault}{\mddefault}{\updefault}$\ddot{\pi}$}}}}
\put(2186,-701){\makebox(0,0)[lb]{\smash{{\SetFigFontNFSS{10}{12.0}{\familydefault}{\mddefault}{\updefault}$\dot{\pi}$}}}}
\put(1206,-1871){\makebox(0,0)[lb]{\smash{{\SetFigFontNFSS{10}{12.0}{\familydefault}{\mddefault}{\updefault}$\pa L+x$}}}}
\put(-1124,-496){\makebox(0,0)[lb]{\smash{{\SetFigFontNFSS{10}{12.0}{\familydefault}{\mddefault}{\updefault}$\pa K$}}}}
\put(-669,-386){\makebox(0,0)[lb]{\smash{{\SetFigFontNFSS{10}{12.0}{\familydefault}{\mddefault}{\updefault}$\pa K'$}}}}
\put(-594,-761){\makebox(0,0)[lb]{\smash{{\SetFigFontNFSS{10}{12.0}{\familydefault}{\mddefault}{\updefault}$a_2$}}}}
\put(-2274,-1671){\makebox(0,0)[lb]{\smash{{\SetFigFontNFSS{10}{12.0}{\familydefault}{\mddefault}{\updefault}$\dot{a_1}$}}}}
\put(-2114,-1986){\makebox(0,0)[lb]{\smash{{\SetFigFontNFSS{10}{12.0}{\familydefault}{\mddefault}{\updefault}$\pa L'+x_0$}}}}
\put(-2624,-741){\makebox(0,0)[lb]{\smash{{\SetFigFontNFSS{10}{12.0}{\familydefault}{\mddefault}{\updefault}$a_1$}}}}
\put(-3294,-1561){\makebox(0,0)[lb]{\smash{{\SetFigFontNFSS{10}{12.0}{\familydefault}{\mddefault}{\updefault}$-b_2+x_0$}}}}
\put(-1199,-1671){\makebox(0,0)[lb]{\smash{{\SetFigFontNFSS{10}{12.0}{\familydefault}{\mddefault}{\updefault}$\dot{a_2}$}}}}
\put(-539,-1951){\makebox(0,0)[lb]{\smash{{\SetFigFontNFSS{10}{12.0}{\familydefault}{\mddefault}{\updefault}$\theta$}}}}
\put(-1734,-641){\makebox(0,0)[lb]{\smash{{\SetFigFontNFSS{10}{12.0}{\familydefault}{\mddefault}{\updefault}$-\dot{b_1}+x_0$}}}}
\put(-2029,-381){\makebox(0,0)[lb]{\smash{{\SetFigFontNFSS{10}{12.0}{\familydefault}{\mddefault}{\updefault}$-\dot{b_2}+x_0$}}}}
\put(-1429,-1871){\makebox(0,0)[lb]{\smash{{\SetFigFontNFSS{10}{12.0}{\familydefault}{\mddefault}{\updefault}$\pa L+x_0$}}}}
\put(-1999,-1231){\makebox(0,0)[lb]{\smash{{\SetFigFontNFSS{10}{12.0}{\familydefault}{\mddefault}{\updefault}$-b_1+x_0$}}}}
\end{picture}%
\end{center}
\caption{$K\cap(L+x)$ (dotted lines) and $K'\cap (L'+x)$ (continuous lines) when $x=x_0$ (left) and when $x=x_0-\ee\te$ (right).}
\label{fig_claim4}
\end{figure}
See Fig.~\ref{fig_claim4}. Let $\pi=\{y\in\Real^2 : (y-a_1)\cdot \te\geq0\}$. The halfplane $\pi$ contains $\Si$ and its boundary contains $a_1$, $a_2$ and $b_1+x_1$. Since $\te=\te'$, we have $b_2+x_1\in\pa\pi$ and $\Om+x_1\subset\pi$. By convexity of the involved polygons,  $K\De K'$ and $(-L+x_1)\De(-L'+x_1)$ are contained in $\Real^2\setminus\pi$. The halfplane $-\pi+x_1+x_0$ has outer normal $\te$ and its boundary contains $\dot{a_1}$, $\dot{a_2}$ and $c_1$. Thus $-\pi+x_1+x_0$ contains $K$, $K'$, $-L+x_1$ and $-L'+x_1$, since its boundary support the four sets at $c_1$.  Summarising, the following inclusions hold for any $x\in\Real^2$:
\begin{equation}\label{inclusioni}\begin{aligned}
&K,K'\subset -\pi+x_1+x_0;& &L+x, L'+x\subset\pi-x_0+x;\\
&K\De K'\subset \Real^2\setminus\pi; & &(L+x)\De (L'+x)\subset\Real^2\setminus(-\pi+x_1+x).
\end{aligned}\end{equation}

Therefore $K\cap(L+x)$ and $K'\cap (L'+x)$ are contained in the strip $N_1(x)=(\pi-x_0+x)\cap(-\pi+x_1+x_0)$, while $K\De K'$ and $(L+x)\De (L'+x)$ do not intersect the strip $N_2(x)=\pi\cap(-\pi+x_1+x)$. Since $N_1(x_0)=N_2(x_0)$ we have $K\cap(L+x_0)=K'\cap(L'+x_0)$. Let $x=x_0-\ee\te$, with $\ee>0$. We have 
\begin{equation}
\big( K\cap(L+x)\big)\Delta  \big(K'\cap(L'+x)\big)\subset N_1(x)\setminus N_2(x)=\dot{\pi}\cup\ddot{\pi},
\end{equation}
where $\dot{\pi}=(\pi-\ee\te)\setminus\pi$ and $\ddot{\pi}=(-\pi+x_1+x_0)\setminus(-\pi+x_1+x_0-\ee\te)$; see Fig.~\ref{fig_claim4}. 

In order to prove $g_{K,L}(x)\neq g_{K',L'}(x)$, we need to distinguish two cases, according to whether $[a_1,a_2]$ is an edge of $K$ (or of $K'$) or not.
Note that  $[a_1,a_2]$ is an edge of $K$  if and only if $[b_1,b_2]$ is an edge of $-L'$. Indeed these conditions are equivalent respectively to $-\te\in N(K,a_1)$ and to $-\te\in N(-L',b_1)$, and  these cones coincide by \eqref{coni_giusti}. 
Similar arguments prove that $[a_1,a_2]$ is an edge of $K'$  if and only if $[b_1,b_2]$ is an edge of $-L$.  We also observe that $[a_1,a_2]$ cannot be an edge of both $K$ and $K'$, since otherwise $\Si=\pa K$, contradicting what has been proved in the lines preceding Claim~\ref{claim1}. 

Assume that $[a_1,a_2]$ is an edge of $K$. In this case $[b_1,b_2]$ is an edge of $-L'$, $\relint [a_1,a_2]\subset K'$  and  $\relint [b_1,b_2]\subset-L$. We have $K\cap\dot{\pi}=\emptyset$, because  $K\subset\pi$. 
When $\ee$ is small, the inclusion \eqref{sit_x0_1} implies 
$(L'+x)\cap \dot{\pi}\subset \inte K'$, by continuity; see Fig.~\ref{fig_claim4}.
Therefore 
\begin{equation*}
\big( \big( K\cap(L+x)\big)\Delta  \big(K'\cap(L'+x)\big)\big)\cap\dot{\pi}=(L'+x)\cap\dot{\pi}.
\end{equation*}
This set is a rectangle of base $\|\dot{b_2}-\dot{b_1}\|$ and height $\ee$, up to  triangles of edge-lengths proportional to $\ee$. Its area is $\ee \|\dot{b_2}-\dot{b_1}\|+o(\ee^2)$. Similar arguments prove that 
\begin{equation*}
\big( \big( K\cap(L+x)\big)\Delta  \big(K'\cap(L'+x)\big)\big)\cap\ddot{\pi}=K\cap(L+x)\cap\ddot{\pi},
\end{equation*}
and that this set has area $\ee \|-b_1+x_0-\dot{a_1}\|+o(\ee^2)$. Therefore we have
\begin{equation*}
g_{K,L}(x)- g_{K',L'}(x)=\ee(\|-b_1+x_0-\dot{a_1}\|-\|\dot{b_2}-\dot{b_1}\|)+o(\ee^2),
\end{equation*}
which, in view of  \eqref{sit_x0_4}, contradicts $g_{K,L}= g_{K',L'}$. 
Assume that $[a_1,a_2]$ is neither an edge of $K$ nor an edge of $K'$ (and, as a consequence, $[b_1,b_2]$ is neither an edge of $-L$ nor an edge of $-L'$). In this case  we have, for $\ee>0$ small, 
\begin{equation*}
\big( \big( K\cap(L+x)\big)\Delta  \big(K'\cap(L'+x)\big)\big)\cap\dot{\pi}=\emptyset,
\end{equation*}
because $(L+x)\cap\dot{\pi}=(L'+x)\cap\dot{\pi}$ (by \eqref{inclusioni}) and both these sets are contained in $K$ and $K'$ (by \eqref{sit_x0_1}). Moreover, \eqref{inclusioni} implies $K\cap\ddot{\pi}=K'\cap\ddot{\pi}$ (because $\ddot{\pi}\subset\pi$ when $\ee$ is small), while \eqref{sit_x0_2} and \eqref{sit_x0_3} imply that $((L+x)\De(L'+x))\cap\ddot{\pi}$ is contained in $B(-b_1+x_0,\de)\cup B(-b_2+x_0,\de)$, for a suitable $\de=\de(\ee)$ positive which tends to $0$ as $\ee$ tends to $0$. If $\ee$ is sufficiently small, $B(-b_2+x_0,\de)$ does not intersect $K\cup K'$ and $B(-b_1+x_0,\de)\cap\ddot{\pi}\subset K, K'$. Thus we have 
\begin{equation*}
\big( \big( K\cap(L+x)\big)\Delta  \big(K'\cap(L'+x)\big)\big)\cap\ddot{\pi}=
\big((L+x)\Delta(L'+x)\big)\cap B(-b_1+x_0,\de)\cap\ddot{\pi}.
\end{equation*}
Arguing as in the last part of the proof of Claim~\ref{claim2} proves $g_{K,L}(x)\neq g_{K',L'}(x)$. We omit the details.

Assume  $\te<\te'$. The formulas~\eqref{intersezione}, \eqref{intersezione_due},  \eqref{intersezionep}, \eqref{intersezione_duep} and \eqref{w1minw2} imply that  either $w_1=\te$ holds or $w_2=\te'$ holds. Assume $w_1=\te$, for instance.  Let $\dot{a_1}$, $\dot{b_1}$, $\dot{a_2}$, $\dot{b_2}$ and $\pi$ be defined as in  case $\te=\te'$. 
Let us prove
\begin{equation}\label{base_piu_lunga2}
\|a_2-a_1\|<\|\dot{a_2}-\dot{a_1}\|.
\end{equation} 
Assume \eqref{base_piu_lunga2} false and define $x_0=a_2+\dot{b_2}$. We have $[a_1,\dot{a_1}]_{\pa K}=[b_1,\dot{b_1}]_{\pa(-L)}+x_1$, $\dot{a_1}=\dot{b_1}+x_1$, $c_1\in [\dot{a_1},\dot{a_2}]_{\pa K}\cap([\dot{b_1},\dot{b_2}]_{\pa (-L)}+x_1)$ and $c_1\neq \dot{a_1}$, because the arguments that prove these relations in the case $\te=\te'$ are valid also in this case. In particular $[\dot{b_1},\dot{b_2}]$ is not a point. The condition $\te<\te'$ implies $-b_2+x_0\notin-\pi+x_1+x_0$.  Since $K$, $K'\subset-\pi+x_1+x_0$, this implies  $-b_2+x_0\notin K\cup K'$. Arguments similar to those used in the case $\te=\te'$  prove that $-\dot{b_1}+x_0$, $-b_1+x_0\notin K\cup K'$. 
Therefore, when $x=x_0-\ee\te$, with $\ee>0$ small, $\big( K\cap(L+x)\big)\Delta  \big(K'\cap(L'+x)\big)$ is contained in a neighbourhood of $a_2$. Arguments similar to those in the last part of the proof of Claim~\ref{claim2} prove $g_{K,L}(x)\neq g_{K',L'}(x)$. We omit the details. This contradiction proves \eqref{base_piu_lunga2}. Arguments similar to those contained in the lines which follow \eqref{base_piu_lunga} prove that \eqref{base_piu_lunga2} implies that $U$ strictly contains an half-circle.

Let us prove $x_1\neq x_2$ arguing by contradiction. If $x_1=x_2$ then $\te=\te'$, by definition. Thus \eqref{rappr_new} and \eqref{rappr_new2} hold and   imply $\Si=\Om+x_1$, contrary to the assumptions of Claim~\ref{claim3}.

Summarising, $U$ may coincide with an half-circle only when $\te=\te'$. When $U$ is an half-circle, $K$ is contained in the strip bounded by the line through $a_1$ orthogonal to $u_1$ and by the line line through $a_2$ orthogonal to $u_1$. Since $[\dot{a_1},\dot{a_2}]$ is contained in this strip and it is parallel to $[a_1,a_2]$, equality holds in \eqref{base_piu_lunga}. Moreover the arcs $[a_1,\dot{a_1}]_{\pa K}$ and $[\dot{a_2},a_2]_{\pa K}$ are line segments contained in the lines bounding the strip. The last part of the claim follows from these observations, \eqref{rappr_new} and \eqref{rappr_new2}.
\end{proof}

\begin{claim}\label{claim4}
The arc $\Si$ is a translate of $\Om$.\end{claim}

\begin{proof}Assume that $\Si$ is not a translate of $\Om$. For $i=1,2$,  let $a_i$, $b_i$ and $u_i$ be as in Claim~\ref{claim1}, let $x_i=b_i-a_i$, let $c_i$ be as in Claim~\ref{claim2} and let $w_i$ be defined as in the lines preceding \eqref{w1minw2}.

If $u\in(u_1,w_1)_{S^1}$, then both \eqref{alt_uno} and \eqref{alt_due} hold, since $K_u$, $K'_u$, $(-L+x_1)_u$ and $(-L'+x_1)_u$ are all contained in  the relative interior of  $\Si\cap (\Om+x_1)$, which is contained in $\pa K\cap\pa K'\cap(-\pa L+x_1)\cap(-\pa L'+x_1)$. 
Let $\widetilde U\subset S^1$ be the maximal arc which contains $(u_1,w_1)_{S^1}$ and such that \eqref{alt_due} holds for each $u\in\widetilde U$. 
Let $\widetilde \Si$  be the  maximal arc of $\pa K\cap \pa (-L'+x_1)$ containing $\cup_{u\in(u_1,w_1)_{S^1}} K_u$ and let  $\widetilde \Om$ be the maximal arc of  $\pa K'\cap \pa (-L+x_1)$ containing $\cup_{u\in(u_1,w_1)_{S^1}} K'_u$.
Since $u_1\neq w_1$, by \eqref{intersezione}, $\widetilde U$ is not a point and neither $\widetilde \Si$ nor $\widetilde \Om$ are points or line segments. Clearly $\widetilde \Si$ and $\widetilde \Om$ contain $a_1$. Moreover, since $\pa K$ and $\pa (-L'+x_1)$ bifurcate at $c_1$ (and the same is true for $\pa K'$ and $\pa (-L+x_1)$), $c_1$ is the upper endpoint of $\widetilde \Si$ and of $\widetilde \Om$ and $\widetilde \Si\neq\pa K$.

If $\widetilde \Si$ is a translate of $\widetilde \Om$, then this translation is the identity, since $\widetilde \Si$ and $\widetilde \Om$ have their upper endpoint in common. On the other hand, \eqref{coni_giusti}, with $j=1$, implies that $\widetilde\Si$ coincides with $\pa K$  and $\widetilde\Om$ coincides with $\pa K'$ in a neighbourhood of $a_1$,. Since $\pa K$ and $\pa K'$ bifurcate at $a_1$,  $\widetilde \Si$ is not a translate of $\widetilde \Om$.

Results analogous to Claim~\ref{claim1} and~\ref{claim3} hold for $\widetilde U$, $\widetilde\Si$ and $\widetilde\Om$. In particular, $\widetilde U$ contains an half-circle.

Let $\bar U\subset S^1$ be the maximal arc which contains $(w_2,u_2)_{S^1}$ and such that \eqref{alt_due} holds for each $u\in\bar U$.
Let $\bar \Si$  be the  maximal arc of $\pa K\cap \pa (-L'+x_2)$ containing $\cup_{u\in(w_2,u_2)_{S^1}} K_u$ and let  $\bar \Om$ be the maximal arc of  $\pa K'\cap \pa (-L+x_2)$ containing $\cup_{u\in(w_2,u_2)_{S^1}} K'_u$.
Also $\bar U$ contains an half-circle. 

For $i=1,2$, let $\widetilde u_i$ and $\bar u_i \in S^1$ be such that $\cl\widetilde U=[\widetilde u_1,\widetilde u_2]_{S^1}$ and $\cl\bar U=[\bar u_1,\bar u_2]_{S^1}$, and let $\widetilde a_1\in\pa K$ and $\bar a_2\in \pa K$ be such that $\widetilde \Si= [\widetilde a_1,c_1]_{\pa K}$ and $\bar \Si=[c_2,\bar a_2]_{\pa K}$.
Let us prove  $\bar u_2\leq\widetilde u_1$. Assume  $\bar u_2>\widetilde u_1$. In this case the sub-arcs $\widetilde \Si$ and $\bar \Si$ of $\pa K$ overlap and contain  the arc $[\widetilde a_1,\bar a_2]_{\pa K}$. The latter is not a point  or a line segment, because arguing as we did in Claim~\ref{claim1} one can prove that $\widetilde U$ contains a line segment orthogonal to $\widetilde u_1$ containing $\widetilde a_1$ and $\bar U$ contains a line segment orthogonal to $\bar u_2$ containing $\bar a_2$. The inclusions
\begin{equation*}
[\widetilde a_1,\bar a_2]_{\pa K}\subset\widetilde \Si\subset\pa(-L'+x_1),\quad\text{and}\quad [\widetilde a_1,\bar a_2]_{\pa K}\subset\bar \Si\subset\pa(-L'+x_2),
\end{equation*}
and the convexity of $\pa(-L')$ imply $x_1=x_2$. This equality contradicts the assumption ``$\Si$ is not a translate of $\Om$'', as shown by Claim~\ref{claim3}, and proves $\bar u_2\leq\widetilde u_1$. Similar arguments prove $\bar u_1\geq \widetilde u_2$. 
These inequalities imply that both $\widetilde U$ and $\bar U$ are half-circles with $\widetilde U\cup\bar U=S^1$.

A description analogous to that of  Claim~\ref{claim3} applies to $\widetilde \Si$ and $\widetilde \Om$ and also to $\bar \Si$ and $\bar \Om$. This description easily implies that  $K$, $K'$, $-L$ and $-L'$ are parallelograms with two edges orthogonal to $\widetilde u_1$ and two edges orthogonal to $v$, for some $v\in S^1$ with $v\neq \widetilde u_1$. It also implies $\len(K_v)=\len((-L')_v)$ and $\len(K'_v)=\len((-L)_v)$. It cannot be $\len(K_{\widetilde u_1})=\len((-L')_{\widetilde u_1})$, because otherwise $\widetilde \Si=\pa K$ and this contradicts what has been proved above. Thus, $\len(K_{\widetilde u_1})=\len(K'_{\widetilde u_1})$ and $\len((-L)_{\widetilde u_1})=\len((-L')_{\widetilde u_1})$, by the synisothesis of $(K,-L)$ and $(K',-L')$.  
Therefore, up to an affine transformation,  $(K,L)$ and $(K',L')$ are trivial associates of $(\cK_3,\cL_3)$ and $(\cK_4,\cL_4)$ (with the defining parameter $m$ equal to $0$), respectively. 
This contradicts the assumptions of Lemma~\ref{archi_simmetrici} and concludes its proof.
\end{proof}

\section{Proof of Theorem~\ref{cov_congiunto_poligoni}}\label{sec_teorema_poligoni}
\begin{proof}
Proposition~\ref{det_curv} implies that  $(K',-L')$ is a pair of polygons synisothetic to $(K,-L)$. In particular, for each $u\in S^1$, either \eqref{alt_uno} or \eqref{alt_due} holds.  
We assume  that $(K,L)$ and $(K',L')$ are not trivial associates and prove that 
\begin{equation}\label{contra}
K=-L+x\quad\text{and}\quad K'=-L'+x',\quad\text{for some $x,x'\in\Real^2$.}
\end{equation}
These identities, together with $K-L=K'-L'$ (which follows by \eqref{support}), imply   $K=K'+(x-x')/2$ and  $L=L'+(x-x')/2$, that is, they prove that $(K,L)$ and $(K',L')$ are trivial associates, concluding the proof. 

In order to prove \eqref{contra}, let $p$ be a vertex of $K$ and $q$ a vertex of $-L$  such that $\relint N(K,p)\cap \relint N(-L,q)\neq\emptyset$. We prove that
\begin{equation}\label{coni_corrispondenti_uguali}
N(K,p)= N(-L,q).
\end{equation}
Let $u_0\in S^1\cap\relint N(K,p)\cap \relint N(-L,q)$ and assume that  \eqref{alt_uno} holds when $u=u_0$. This condition implies that there exist $y$, $y'\in\Real^2$ such that  $K$ and $K'+y$ coincide in a neighbourhood of $p$, while $-L$ and $-(L'+y')$ coincide in a neighbourhood of $q$. Formulas~\eqref{support} and~\eqref{facce_corpodiff} imply
\[
p+q=K_{u_0}+(-L)_{u_0}=K'_{u_0}+(-L')_{u_0}=(p-y)+(q+y'),
\]
that is $y=y'$.
We apply Lemma~\ref{archi_simmetrici} to $(K,L)$ and $(K'+y,L'+y)$, with $U$ chosen so that it contains $N(K,p)\cap N(-L,q)$.
If $\Si$ and $\Om$ are defined as in the
statement of Lemma~\ref{archi_simmetrici}, then they are not points nor line segments.
This lemma implies that
$\Si$ is a translate  of $\Om$,
which yields \eqref{coni_corrispondenti_uguali}. Similar arguments prove \eqref{coni_corrispondenti_uguali} when \eqref{alt_due} replaces \eqref{alt_uno}. In this case we apply Lemma~\ref{archi_simmetrici} to $(K,L)$ and $(-L'+y,-K'+y)$, where   $y\in\Real^2$ is chosen so that $K$ and $-L'+y$ coincide in a neighbourhood of $p$. 

What has been proved so far implies that to each edge $E$ of $K$ it corresponds an edge $F$ of $-L$ with equal outer normal, and vice versa. To prove that $K$ is a translate of $-L$ it suffices to show that 
\begin{equation}\label{lati_lunghi_uguali}
\len(E)=\len(F).
\end{equation}
Let $E=[x_1, x_2]$ and $F=[y_1,y_2]$. We may label the vertices in such a way that $\inte N(K,x_i)\cap \inte N(-L,y_i)\neq\emptyset$, for each $i=1,2$.
Let $u_0$ be the unit outer normal to $K$ at $E$ and assume that \eqref{alt_uno} holds when $u=u_0$. One proves, arguing as above, that there exists $y\in\Real^2$  such that  $E$ is an edge of $K'+y$ with outer normal $u_0$ and $F$ is an edge of $-(L'+y)$ with outer normal $u_0$. We apply Lemma~\ref{archi_simmetrici} to $(K,L)$ and $(K'+y,L'+y)$, with $U$ chosen so that it contains $u_0$.
What has been proved above implies that  $K$ and $K'+y$ coincide in a neighbourhood of $E$, that $-L$ and $-(L'+y)$ coincide in a neighbourhood of $F$, and that $U$ is not a point. 
 If $\Si$ and $\Om$ are defined as in the
statement of Lemma~\ref{archi_simmetrici}, then they are not points nor line segments.
This lemma  implies that $\Si$ is a translate of $\Om$,
which yields \eqref{lati_lunghi_uguali}. Similar arguments get the same conclusion  when \eqref{alt_due} substitutes \eqref{alt_uno},  and similar arguments also prove   that $K'$ is a translate of $-L'$. 
\end{proof}

\begin{proof}[Proof of Corollary~\ref{symmetry_crosscov}]
First we prove the corollary assuming $z=0$.
Assume $g_{K,L}(x)=g_{K,L}(-x)$ for each $x\in\Real^2$.  This is equivalent to  $g_{K,L}(x)=g_{L,K}(x)$ for each $x\in\Real^2$, since $g_{L,K}(x)=g_{K,L}(-x)$. 
We claim that  there exist no affine transformation $\cT$ and no different  indices $i,j$, with either $i,j\in\{1,2\}$ or  $i,j\in\{3,4\}$, such that $(\cT K,\cT L)$ and $(\cT L,\cT K)$ are trivial associates of $(\cK_i,\cL_i)$ and $(\cK_j,\cL_j)$, respectively. Indeed, if this claim is false, then $(\cK_i,\cL_i)$ is a trivial associate of  $(\cL_j,\cK_j)$, because being trivial associates is a transitive property. However, when $i\neq j$,  $(\cK_i,\cL_i)$ is not a trivial associate of  $(\cL_j,\cK_j)$, because $\cK_i$ is not a translate of $-\cK_j$ or of $\cL_j$. 
This claim and Theorem~\ref{cov_congiunto_poligoni} imply that $(K,L)$ is a trivial associate of $(L,K)$. It is immediate to understand that this happens exactly when  $K=-K+y$ and $L=-L+y$, for some $y\in\Real^2$, (that is, $y/2$ is the center of  $K$ and of $L$) or when $K=L$.
The converse implication follows from the identities $g_{K,L}(x)=g_{-K+y,-L+y}(-x)=g_{L,K}(-x)$, valid for any $x,y\in\Real^2$.

The proof for $z\neq0$ follows from the one for $z=0$ applied to $g_{K,L+z}$, since  $g_{K,L}(z+x)= g_{K,L+z}(x)$ and $g_{K,L}(z-x)= g_{K,L+z}(-x)$.
\end{proof}

\textsc{Acknowledgements.}
We are extremely grateful to G.~Averkov and  R.~J.~Gardner for reading large parts of this manuscript and suggesting many arguments that simplified and clarified some proofs.
We also thank  P.~Mani-Levitska for giving us his unpublished note~\cite{ML}.

\bibliographystyle{amsplain}

\end{document}